\documentclass[leqno,twoside,11pt]{amsart}
\usepackage{amsmath,amssymb,latexsym,cite,mathrsfs,amsfonts}
\usepackage{enumitem,graphicx}
\usepackage[hidelinks,linktocpage,pdfpagelabels,
bookmarksnumbered,bookmarksopen]{hyperref}
\usepackage[english]{babel}
\usepackage{comment}
\usepackage{tikz}

\usepackage[left=2.9cm,right=2.9cm,top=2.8cm,bottom=2.8cm]{geometry}

\usepackage{verbatim} 

\numberwithin{equation}{section}

\newtheorem{theorem}{Theorem}[section]
\newtheorem{proposition}[theorem]{Proposition}
\newtheorem{lemma}[theorem]{Lemma}
\newtheorem{corollary}[theorem]{Corollary}
\newtheorem{remark}[theorem]{Remark}

\theoremstyle{definition}
\newtheorem{definition}[theorem]{Definition}
\theoremstyle{plain}

\newcommand\restr[2]{{
  \left.\kern-\nulldelimiterspace 
  #1 
  \vphantom{\big|} 
  \right|_{#2} 
  }}

\title[Mountain pass levels on cones]{Nehari type characterizations of the mountain pass level on cones}

\author[D. S. Pereira]{Denilson S. Pereira}
\address[Denilson S. Pereira]{Unidade Acadêmica de Matemática, Universidade Federal de Campina Grande, Rua Aprígio Veloso 882, 58429-900 Campina Grande, PB, Brazil}
\email{denilsonsp@mat.ufcg.edu.br}
\author[H. Ramos Quoirin]{Humberto Ramos Quoirin}
\address[Humberto Ramos Quoirin]{CIEM (CONICET), FaMAF, Universidad Nacional de Córdoba, Medina Allende s/n, Ciudad Universitaria, 5000 Córdoba, Argentina}
\email{humbertorq@gmail.com}
\author[J. R. Santos Júnior]{João R. Santos Júnior\textsuperscript{*}}
\address[João R. Santos Júnior]{Faculdade de Matemática, Instituto de Ciências Exatas e Naturais, Universidade Federal do Pará, Rua Augusto Corrêa 1, Guamá, 66075-110 Belém, PA, Brazil}
\email{joaojunior@ufpa.br}
\author[K. Silva]{Kaye Silva\textsuperscript{\ensuremath{\dagger}}}
\address[Kaye Silva]{Instituto de Matemática e Estatística, Universidade Federal de Goiás, Rua Samambaia, 74001-970 Goiânia, GO, Brazil}
\email{kayesilva@ufg.br}
\thanks{\textsuperscript{*}João R. Santos Júnior is the corresponding author (\texttt{joaojunior@ufpa.br}) and holds a CNPq research productivity fellowship (grant 304340/2025-1).}
\thanks{\textsuperscript{\ensuremath{\dagger}}Kaye Silva holds a CNPq research productivity fellowship (grant 300849/2025-7).}

\begin{document}

\begin{abstract} We relate the mountain pass level of a functional to its infimum over the part of the Nehari set formed by local maximizers of the fiber maps. {Unlike in the classical setting where every fiber map has a unique global maximizer,} the fiber maps need not all have a maximizer: working on each path-connected component of the cone of directions where they do, we show that the two levels coincide, for merely continuous functionals and without any compactness assumption. We also introduce a notion of mountain pass geometry with respect to a cone, which may hold when the classical one fails. {This geometry includes the zero-level case when a suitable separating set exists.} Under the Palais--Smale condition, each component carries a critical point at its own mountain pass level. Applications include indefinite and concave-convex problems with sign-changing weights, the affine $p$-Laplacian, nonlocal problems, and a problem related to the Fu\v{c}\'{\i}k spectrum.

\end{abstract}

\maketitle

\bigskip

\begin{center}
	\begin{minipage}{12cm}
		\tableofcontents
	\end{minipage}
\end{center}

\section{Introduction}

\begin{figure}[ht]
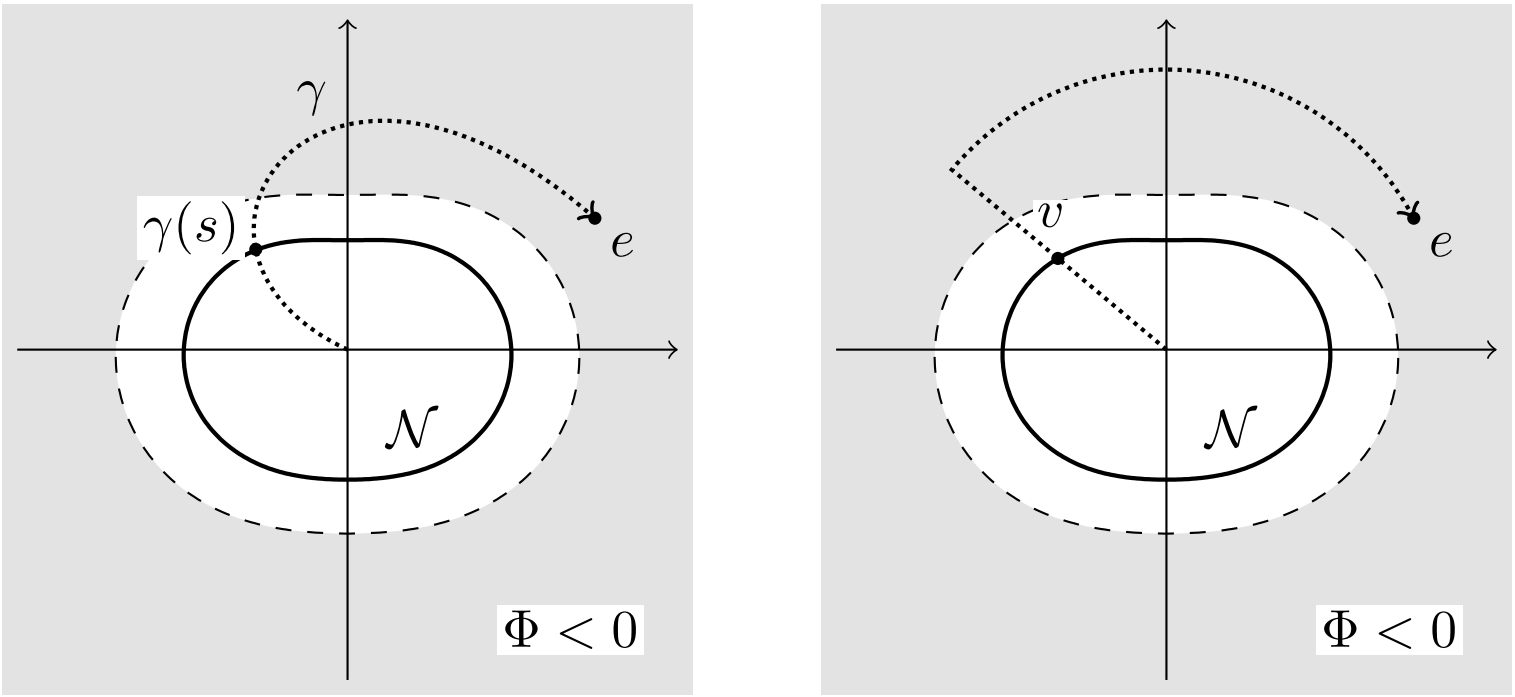

\centering

\caption{The situation behind \eqref{equali}, drawn for $\Phi=\frac{1}{2}g(\theta)r^2-\frac{1}{4}r^4$ in polar coordinates, with $g>0$, so that every fiber map has a unique critical point, which is a global maximizer. The Nehari set $\mathcal{N}$ (solid) is a closed curve around the origin, the curve $\{\Phi=0\}$ is dashed, and the region $\{\Phi<0\}$ is shaded. Left: every path $\gamma$ from $0$ to $e$ crosses $\mathcal{N}$, which gives $c(e)\ge\inf_{\mathcal{N}}\Phi$. Right: the path that follows the ray through $v\in\mathcal{N}$ and then stays in $\{\Phi\le0\}$ has maximum $\Phi(v)$, which gives $c(e)\le\inf_{\mathcal{N}}\Phi$.}
\label{figclassic}
\end{figure}

Let $X$ be a real normed space and let $\Phi\in C^1(X)$ be such that $\Phi(0)=0$. A classical way of finding nontrivial critical points of $\Phi$ is to study its \emph{fiber maps}
\begin{equation*}
	\varphi_u(t)=\Phi(tu), \qquad t>0, \quad u\in X\setminus\{0\}.
\end{equation*}
A point $tu$ belongs to the Nehari set $\mathcal{N}=\{u\in X\setminus\{0\}: \Phi'(u)u=0\}$ if, and only if, $t$ is a critical point of $\varphi_u$. Thus $\mathcal{N}$ contains all nontrivial critical points of $\Phi$, and its structure is determined by the shape of the fiber maps \cite{SW}. {In the classical setting where every fiber map has a unique critical point, which is a global maximizer,} the minimization of $\Phi$ over $\mathcal{N}$ is closely related to the mountain pass theorem \cite{AR}. Indeed, given $e\in X$, let
\begin{equation*}
	c(e)=\inf_{\gamma\in\Gamma(e)}\max_{t\in[0,1]}\Phi(\gamma(t)), \qquad \Gamma(e)=\{\gamma\in C([0,1],X): \gamma(0)=0,\ \gamma(1)=e\}.
\end{equation*}
Under suitable conditions, which include $\dim X\ge2$ (see Remark \ref{dimone}), one has
\begin{equation}\label{equali}
	c(e)=\inf_{\mathcal{N}}\Phi=\inf_{u\in X\setminus\{0\}}\max_{t>0}\Phi(tu)
\end{equation}
for every $e\neq0$ such that $\Phi(e)\le0$. The proof of \eqref{equali} relies on the structure of $\mathcal{N}$. Denoting by $t(u)$ the critical point of $\varphi_u$, we have $\mathcal{N}=\{t(u)u: u\in X\setminus\{0\}\}$ and $\max_{t>0}\Phi(tu)=\Phi(t(u)u)$, which gives the second equality. The inequality $c(e)\ge\inf_{\mathcal{N}}\Phi$ holds because every path $\gamma\in\Gamma(e)$ crosses $\mathcal{N}$ (Figure \ref{figclassic}, left): the map $s\mapsto t(\gamma(s))$ is larger than $1$ shortly after the last time at which $\gamma$ vanishes, since $\mathcal{N}$ stays away from the origin, and it is smaller than $1$ at $s=1$, since $e$ lies beyond $\mathcal{N}$ on its ray; the intermediate value theorem then yields some $s$ with $\gamma(s)\in\mathcal{N}$. Conversely, given $v\in\mathcal{N}$, one joins $0$ to $e$ by following the ray through $v$ until $\Phi$ vanishes, and then moving inside the set of points lying beyond $\mathcal{N}$ on their ray at which $\Phi\le0$ (Figure \ref{figclassic}, right). This set is path-connected if $\dim X\ge2$, and the maximum of $\Phi$ along the resulting path is $\Phi(v)$, which yields $c(e)\le\inf_{\mathcal{N}}\Phi$. To the best of our knowledge, \eqref{equali} goes back to Ding and Ni \cite{DiNi}; see also Ni and Takagi \cite{NiTa}, Rabinowitz \cite{R3} and Willem \cite[Theorem 4.2]{Wi}. It shows that the least energy level on the Nehari set is a mountain pass level, which in particular does not depend on $e$.

In many problems, however, the fiber maps do not behave uniformly: $\varphi_u$ has a local maximizer only for $u$ in an open cone
\begin{equation*}
	Y=\{u\in X\setminus\{0\}: \varphi_u \mbox{ has a local maximizer}\},
\end{equation*}
which may be a proper subset of $X\setminus\{0\}$ and may be disconnected. This is the case for indefinite problems such as $-\Delta u=\lambda u+b(x)|u|^{p-2}u$ with a sign-changing weight $b$, where, for suitable values of $\lambda$, $\varphi_u$ has a maximum if, and only if, $\int_\Omega b|u|^p>0$ (Subsection \ref{subind}). The same occurs for concave-convex problems with sign-changing weights (Subsection \ref{subcc}) and for jumping nonlinearities (Subsection \ref{subfucik}); see Figure \ref{fig1}, to be compared with Figure \ref{figclassic}. For such problems one usually minimizes $\Phi$ over the set $\mathcal{N}^-$ of points $t^-(u)u$, where $u\in Y$ and $t^-(u)$ is the local maximizer of $\varphi_u$ \cite{BW,BZ,FRS,Ta}. The relation between this minimization and minimax levels is then much less clear, since both steps of the classical proof of \eqref{equali} break down. First, the intermediate value argument requires the path to stay in $Y$, where $t^-$ is defined, but a path from $0$ to $e$ may leave $Y$, so that it need not cross $\mathcal{N}^-$. Second, the set of points of $Y$ lying beyond $\mathcal{N}^-$ at which $\Phi\le0$ may be disconnected, so that the paths used to bound $c(e)$ from above cannot be built in general.

In this paper we show that these difficulties disappear if one works on each path-connected component $W$ of $Y$. Our main contributions are the following.
\begin{enumerate}[label=(\roman*)]
	\item \emph{A minimax characterization on each component} (Theorem \ref{thmmain}). Assuming only that, for $u\in W$, the fiber map $\varphi_u$ has a unique local maximizer and a negative limit at infinity, and that $\inf_{\mathcal{N}^-\cap W}\Phi\ge0$,
	\begin{equation*}
		c(e)=\inf_{\mathcal{N}^-\cap W}\Phi=\inf_{u\in W}\max_{t>0}\Phi(tu)
	\end{equation*}
	for every $e\in W$ lying beyond $\mathcal{N}^-\cap W$ and such that $\Phi(e)\le0$. Here $\Phi$ is merely continuous, no compactness is assumed, no regularity of $\mathcal{N}^-\cap W$ is required, no condition is imposed at the boundary of $W$, and the paths in $\Gamma(e)$ are not required to stay in $W$, nor in $Y$. Each component of $Y$ thus carries its own mountain pass level, and these levels may differ from one component to another (Remark \ref{remg}).
{
	\item \emph{Mountain pass geometry with respect to a cone} (Definition \ref{defMPG} and Theorem \ref{thmMPG}). In this notion the separating set is required to lie in the cone, while the paths are not. Under a condition on the fiber maps near $\partial W$, the set $\mathcal{N}^-\cap W$ is a closed separating set on which the infimum of $\Phi$ is the mountain pass level. This geometry may hold even when the classical one fails: the planar example of Theorem \ref{thm3.1} has no local minimizers (Remark \ref{MPGC} and Figure \ref{fig1}), while the nonlocal elliptic example of Theorem \ref{thm6} establishes failure of the classical mountain pass geometry without asserting the absence of local minimizers. We also explain why this phenomenon requires the set of directions where the nonlinear term is positive to be disconnected (Remark \ref{local}). {The geometry also accommodates a vanishing mountain pass level, under an additional radial separation condition.}

}
	\item \emph{Critical points with a minimax characterization} (Subsection \ref{subcrit}). When $\Phi$ is of class $C^1$, the level $c_W^-$ is a critical value under the Palais--Smale condition (Corollary \ref{existence}), and minimizers of $\Phi$ on $\mathcal{N}^-\cap W$ are critical points, of Morse index one if $\Phi$ is of class $C^2$ (Proposition \ref{natural}). In particular, if $Y$ has several components, each of them carries a critical point, which is a mountain pass point relative to that component (Corollary \ref{multiplicity}). Critical points obtained from different cones already appear in \cite{FRS}, and in problems where nonnegative and nonpositive solutions are sought separately; what is new here is their characterization as mountain pass points relative to the corresponding components.
	\item \emph{Applications} (Section \ref{secapp}). We obtain a minimax characterization of the Nehari levels of indefinite and concave-convex problems with sign-changing weights, including an extremal case where the mountain pass level vanishes {and the cone mountain pass geometry still holds}, an existence result for the affine $p$-Laplacian with an indefinite weight, whose energy is not coercive, and the elliptic problem mentioned in (ii).
\end{enumerate}

\subsection{Setting and main result}
Our main results do not require $\Phi$ to be differentiable: in Theorems \ref{thmmain} and \ref{thmMPG}, $X$ is a real normed space and $\Phi\in C(X)$ satisfies $\Phi(0)=0$. Differentiability is only needed to obtain critical points (Subsection \ref{subcrit}). Let $Y$ be as above, and assume that $Y$ is nonempty and open; it is then an open cone. Let $W$ be a path-connected component of $Y$. If $u\in Y$ and $\varphi_u$ has only one local maximizer, we denote it by $t^-(u)$. This is the case for every $u\in W$ under (H1') below, and for every $u\in Y$ in all our applications. We set
\begin{equation*}
	\mathcal{N}^-=\{t^-(u)u: u\in Y,\ t^-(u)\mbox{ is defined}\}, \qquad c_W^-=\inf_{\mathcal{N}^-\cap W}\Phi.
\end{equation*}
Given $e\in W$, we write $c_W(e)=c(e)$ and $\Gamma=\Gamma(e)$; the subscript $W$ only records that $e\in W$. Our first main result requires the following conditions, which only concern the fiber maps of directions in $W$:
\begin{enumerate}
	\item[(H1')] For any $u\in W$ the map $\varphi_u$ has exactly one local maximizer.
	\item[(H2')] $\displaystyle\lim_{t\to\infty}\varphi_u(t)<0$ for every $u\in W$.
	\item[(H3')] $c_W^-\ge 0$.
\end{enumerate}
Under (H1'), each $\varphi_u$ with $u\in W$ has at most one local minimizer on each side of $t^-(u)$, and the map $u\mapsto t^-(u)$ is automatically continuous in $W$ (Lemma \ref{tcont}), so that no regularity of $\mathcal{N}^-\cap W$ has to be assumed. By (H2'), for every $u\in W$ there exists $t>t^-(u)$ such that $\Phi(tu)<0$, so there are points $e\in W$ such that $\Phi(e)\le 0$ and $t^-(e)<1$, i.e. lying beyond $\mathcal{N}^-\cap W$ on their ray.

\begin{theorem}\label{thmmain}
Let $W$ be a path-connected component of $Y$ and let $e\in W$ be such that $\Phi(e)\le 0$ and $t^-(e)<1$. Assume (H1')--(H3'). Then
\begin{equation*}
	c_W(e)=c^-_W=\inf_{u \in W} \max_{t>0} \Phi(tu).
\end{equation*}
In particular, $c_W(e)$ does not depend on the choice of $e\in W$.
\end{theorem}

Theorem \ref{thmmain} extends \eqref{equali} to fiber maps that do not behave uniformly. When $Y=X\setminus\{0\}$ and $\dim X\ge 2$, it yields \eqref{equali} with $\mathcal{N}$ replaced by $\mathcal{N}^-$ (Corollary \ref{cor1}). No condition on the behavior of the fiber maps near $\partial W$ is needed. Such a condition will only be used to ensure that $\mathcal{N}^-\cap W$ itself separates $0$ from $e$ (Theorem \ref{thmMPG}). Without it, a path from $0$ to $e$ may reach the region beyond $\mathcal{N}^-\cap W$ through the cone over $\partial W$ without meeting $\mathcal{N}^-\cap W$; it then reaches the level $c_W^-$ on that cone (Proposition \ref{lower}). Theorem \ref{thmmain} also covers the degenerate case $c_W^-=0$, where the mountain pass level coincides with $\Phi(0)$ (Remark \ref{zerolevel}). Condition (H3') cannot be dropped: since every path in $\Gamma$ starts at $0$, we always have $c_W(e)\ge\Phi(0)=0$.

\subsection{Mountain pass geometry with respect to a cone}
Let us now make precise what we mean by mountain pass geometry with respect to a cone.

\begin{definition}\label{defMPG}
Let $W\subset X\setminus\{0\}$ be an open cone and let $e\in W$ be such that $\Phi(e)\le 0$. We say that $\Phi$ has the \emph{mountain pass geometry with respect to $W$} (between $0$ and $e$) if there exists a set $F\subset W$ such that:
\begin{enumerate}[label=(\alph*)]
	\item every $\gamma\in C([0,1],X)$ with $\gamma(0)=0$ and $\gamma(1)=e$ intersects $F$;
	\item $\displaystyle\inf_F\Phi{\ge0}$.
\end{enumerate}
{If $\inf_F\Phi=0$, we additionally require $F$ to be closed in $X$ and $e\notin F$. We call this the \emph{zero-level} case; otherwise the geometry is of \emph{positive level}.}
\end{definition}

This extends the classical mountain pass geometry \cite{AR}: if $\Phi\ge\alpha>0$ on a sphere $\partial B_\rho$ with $\rho<\|e\|$, then (a) and (b) hold with $W=X\setminus\{0\}$ and $F=\partial B_\rho$. As in Ghoussoub and Preiss \cite{GhPr}, the separating set need not be a sphere. The point of Definition \ref{defMPG} is that $F$ is required to lie inside the cone $W$, while the paths are allowed to leave $W$. {In the zero-level case, closedness and the exclusion of $e$ ensure genuine separation of the two endpoints, without requiring $F$ to be homeomorphic to a sphere.}

In order to show that $\mathcal{N}^-\cap W$ is a separating set, we need some control of the fiber maps of directions close to $\partial W$. Let $S=\{u\in X: \|u\|=1\}$ be the unit sphere of $X$. We shall assume either
\begin{enumerate}
	\item[(H4')] If $(u_n)\subset S\cap W$ and $u_n\to u\in\partial W$, then $t^-(u_n)\to\infty$,
\end{enumerate}
which means that the points of $\mathcal{N}^-\cap W$ go to infinity as their direction approaches $\partial W$, or the following strengthening of (H1'), which holds in many applications, where $\varphi_u$ has no local minimizer after $t^-(u)$:
\begin{enumerate}
	\item[(H1'')] Condition (H1') holds and, for any $u\in W$, the map $\varphi_u$ has no local minimizer in the interval $(t^-(u),\infty)$.
\end{enumerate}
When $c_W^->0$, condition (H1'') implies (H4') (Lemma \ref{infinity}).
{Only when $c_W^-=0$, we shall impose the additional condition
\begin{enumerate}
\item[(H5')] $\displaystyle\inf_{u\in S\cap W}t^-(u)>0$.
\end{enumerate}
It is equivalent to $\operatorname{dist}(0,\mathcal{N}^-\cap W)>0$ and prevents the Nehari set from approaching the origin. For $c_W^->0$ this property already follows from the continuity of $\Phi$ at $0$. Under (H5'), (H1'') also implies (H4') at the zero level (Lemma \ref{infinity}).}

\begin{theorem}\label{thmMPG}
Let $W$ be a path-connected component of $Y$ and let $e\in W$ be such that $\Phi(e)\le 0$ and $t^-(e)<1$. Assume that {either $c_W^->0$, or $c_W^-=0$ and (H5') holds,} and that either (H1'') holds, or (H1') and (H4') hold. Then $\mathcal{N}^-\cap W$ is closed in $X$, and every $\gamma\in C([0,1],X)$ with $\gamma(0)=0$ and $\gamma(1)=e$ intersects $\mathcal{N}^-\cap W$. In particular, $\Phi$ has the mountain pass geometry with respect to $W$, with $F=\mathcal{N}^-\cap W$. If, moreover, (H2') holds, then
\begin{equation*}
	\inf_F\Phi=c_W(e).
\end{equation*}
\end{theorem}

The last assertion of Theorem \ref{thmMPG} says that the infimum of $\Phi$ on the closed separating set $F$ equals the mountain pass level. This is the limiting case of the mountain pass principle of Ghoussoub and Preiss \cite{GhPr}, which we use in Corollary \ref{existence} to show that $c_W^-$ is a critical value under the Palais--Smale condition. {The zero-level alternative requires the extra radial bound to make the separating set closed and prevent convergence to the trivial critical point.}

\subsection{Related works}
The decomposition of the Nehari set according to the nature of the critical points of the fiber maps, together with the minimization of $\Phi$ over $\mathcal{N}^-$, has been widely used for problems whose fiber maps do not behave uniformly, such as nonhomogeneous problems \cite{Ta}, indefinite problems with sign-changing weights \cite{BZ}, and concave-convex problems \cite{BW}, which go back to \cite{ABC}; for indefinite problems we also refer to \cite{AlTa,BCN}. The extremal parameter values up to which this approach applies were characterized in \cite{I1} by means of nonlinear Rayleigh quotients, and we refer to \cite{SW} for a survey of the Nehari method. In \cite{Ta}, minimization over $\mathcal{N}^-$ is combined with the mountain pass principle of Ghoussoub and Preiss \cite{GhPr}. An abstract treatment of the splitting of the Nehari set into its parts made of local minimizers and local maximizers of the fiber maps, when all the fiber maps have two critical points, is given in \cite{LRS2} for functionals of class $C^1$ away from the origin (see also \cite{LRS1}).

{Generalizations of the classical mountain pass theorem were developed by Pucci and Serrin {\cite{PuSe84,PuSe85}}, including a zero-altitude formulation; they also studied the structure of the resulting critical set {\cite{PuSe87}}. Ghoussoub and Preiss {\cite{GhPr}} established a closed-set separation principle which also covers the zero-altitude case.}

The closest work to ours is \cite{FRS}, where the Nehari method is restricted to an open cone $Y$ on which the fiber maps behave uniformly. There, ground states relative to $Y$ are obtained by minimizing $\Phi$ over $\mathcal{N}\cap Y$, for several geometries of the fiber maps, under weak lower semicontinuity and compactness assumptions, and several applications to equations and systems are given. The condition labeled (H1) in \cite{FRS}, which requires $\varphi_u$ to have a unique critical point that is a global maximizer, is a particular case of (H1''). The present paper addresses a different question: Theorem \ref{thmmain} identifies $\inf_{\mathcal{N}^-\cap W}\Phi$ with a minimax level, for merely continuous functionals and without any compactness; it works on each path-connected component of $Y$ rather than on the whole cone; and Definition \ref{defMPG} and Theorem \ref{thmMPG} make precise the notion of mountain pass geometry with respect to a cone, a term used informally in \cite{FRS} to describe fiber maps having a maximum. Variational structures over cones were also considered in \cite{DL}, in the framework of linking theorems.

\subsection{Applications and organization of the paper}
The applications in Section \ref{secapp} illustrate the different features of our results:
\begin{itemize}
	\item \emph{Planar examples} (Subsection \ref{subfin}): several cones of mountain pass type, at different levels (Theorem \ref{thm3} and Remark \ref{remg}) or at the same level (Theorem \ref{thmsat}), and a functional without local minimizers that has the mountain pass geometry with respect to a cone but not the classical one (Theorem \ref{thm3.1} and Figure \ref{fig1}).
	\item \emph{Indefinite and concave-convex problems} (Subsections \ref{subind} and \ref{subcc}): a minimax characterization, on each component of $Y$, of Nehari levels whose existence is known \cite{FRS,RS,RSiS}, including an extremal case where the mountain pass level vanishes (Remark \ref{zerolevel}; see also \cite{RSS}) {and the cone mountain pass geometry holds at zero level}.
	\item \emph{The affine $p$-Laplacian with an indefinite weight} (Subsection \ref{subaff}), whose energy is not coercive and is only known to be of class $C^1$ away from the origin.
	\item \emph{A nonlocal problem} (Subsection \ref{subnl}) for which the classical mountain pass geometry fails while the one with respect to a cone holds, together with an obstruction to this phenomenon for local power-type nonlinearities (Remark \ref{local}).
	\item \emph{Disconnected cones} (Subsections \ref{subfr} and \ref{subfucik}): finite-rank problems whose cone has $2^m$ components, and a non-even local problem related to the Fu\v{c}\'{\i}k spectrum, whose cone is disconnected. In these examples the critical points can also be obtained by other means (explicitly, or by considering nonnegative and nonpositive solutions separately \cite{DHI}); they illustrate how Corollary \ref{multiplicity} characterizes each of them as a mountain pass point relative to its component. For the latter problem we show that the first nontrivial curve of the Fu\v{c}\'{\i}k spectrum is exactly the threshold below which the cone is disconnected and the mountain pass level depends on the endpoint (Theorem \ref{thmfucik} and Remark \ref{remfucik}).
	\item \emph{Prescribed scalar curvature} on manifolds of negative Yamabe type (Subsection \ref{subcurv}), where Theorem \ref{thmmain} holds at the critical Sobolev exponent, without any compactness.
\end{itemize}
Section \ref{secproofs} contains the proofs of Theorems \ref{thmmain} and \ref{thmMPG}, the results on critical points and multiplicity, and the classical setting, including the role of the assumption $\dim X\ge2$ and a translation analogue of our results for Neumann problems (Remark \ref{transl}). In Section \ref{secapp}, each subsection starts with the setting and the statement of the corresponding result, followed by its proof.

\section{Proofs and further results}\label{secproofs}

In this section we prove Theorems \ref{thmmain} and \ref{thmMPG}, and we establish further results on critical points and multiplicity (Subsection \ref{subcrit}) and on the classical setting (Subsection \ref{subclass}). Throughout, $W$ is a path-connected component of $Y$ and $e\in W$ satisfies $\Phi(e)\le 0$ and $t^-(e)<1$.

\subsection{Preliminary results}

\begin{lemma}\label{pathcom} Let $W\subset Y$ be a path-connected component of $Y$. Then $W$ is an open path-connected cone.
\end{lemma}
\begin{proof} A ball centered at $w\in W$ and contained in the open set $Y$ is path-connected, hence contained in $W$. For $t>0$, the segment from $w$ to $tw$ lies in the ray of $w$, hence in the cone $Y$, so $tw\in W$.
\end{proof}

Since $\varphi_{\sigma u}(t)=\varphi_u(\sigma t)$ for $\sigma,t>0$, we have $t^-(\sigma u)=\sigma^{-1}t^-(u)$ for every $u\in Y$ and $\sigma>0$. In particular, since $W$ is a cone, $\mathcal{N}^-\cap W=\{t^-(u)u: u\in W\}$. Moreover, $S\cap W$ is path-connected: if $\eta$ is a path in $W$ joining $u,v\in S\cap W$, then $\eta/\|\eta\|$ is a path in $S\cap W$ joining $u$ and $v$. Finally, $\partial W\cap Y=\emptyset$, since a ball contained in $Y$ and centered at a point of $\partial W$ would meet $W$, and hence be contained in $W$. Thus, for $v\in\partial W\setminus\{0\}$, the map $\varphi_v$ has no local maximizer in $(0,\infty)$, and therefore
\begin{equation}\label{nomax}
	\varphi_v(s)\ge\varphi_v(\tau) \quad \mbox{whenever } 0<\tau\le s \mbox{ and } \varphi_v(\tau)\ge0, \qquad v\in\partial W\setminus\{0\}.
\end{equation}
Indeed, otherwise the maximum of $\varphi_v$ on $[0,s]$, where $\varphi_v(0):=\Phi(0)=0$, would be larger than $\varphi_v(s)$ and at least $\varphi_v(\tau)\ge0$, so it would be achieved at a point of $(0,s)$ (at $\tau$, if it vanishes), i.e. at a local maximizer of $\varphi_v$.

\begin{lemma}\label{tcont}
Assume (H1') and let $u\in W$. Then $\varphi_u$ has at most one local minimizer in $(0,t^-(u))$, at most one in $(t^-(u),\infty)$, and no other local extrema besides $t^-(u)$. Moreover, the map $u\mapsto t^-(u)$ is continuous in $W$.
\end{lemma}

\begin{proof}
By (H1'), every local extremum of $\varphi_u$ other than $t^-(u)$ is a local minimizer. If $s_1<s_2$ were local minimizers on the same side of $t^-(u)$, a maximum point $s$ of $\varphi_u$ on $[s_1,s_2]$ would be a local maximizer different from $t^-(u)$: this is clear if $s\in(s_1,s_2)$, and if, say, $s=s_1$, then $\varphi_u$ is constant on some $[s_1,s_1+\delta]$, whose interior points are local maximizers. This proves the first assertion.

Recall that a continuous function without local extrema in an open interval is strictly monotone there. Let $t_0=t^-(u)$ and let $\varepsilon\in(0,t_0)$ be such that $t_0$ is the only local extremum of $\varphi_u$ in $[t_0-\varepsilon,t_0+\varepsilon]$ (there are at most three), so that $\varphi_u(t_0\pm\varepsilon)<\varphi_u(t_0)$. By continuity, there exists $\delta>0$ such that $v\in W$ and $\Phi((t_0\pm\varepsilon)v)<\Phi(t_0v)$ whenever $\|v-u\|<\delta$. For such $v$, $\varphi_v$ has a local maximizer in $(t_0-\varepsilon,t_0+\varepsilon)$, which equals $t^-(v)$ by (H1').
\end{proof}

The continuity of $t^-$ is usually derived from the implicit function theorem, which requires some differentiability of $\Phi$ and the nondegeneracy condition $\varphi_u''(t^-(u))<0$. Lemma \ref{tcont} shows that, under (H1'), it holds for merely continuous functionals.

\begin{lemma}\label{shape}
Assume (H1'), (H2') and (H3'), and let $u\in W$. Then $t^-(u)$ is a global maximizer of $\varphi_u$ in $(0,\infty)$. Moreover, if $r_1(u)$ denotes the smallest $r\ge t^-(u)$ such that $\varphi_u(r)=0$, then $r_1(u)$ is the only zero of $\varphi_u$ in $[t^-(u),\infty)$ and $\varphi_u<0$ in $(r_1(u),\infty)$.
\end{lemma}

\begin{proof}
By Lemma \ref{tcont}, in $(0,t^-(u))$ the map $\varphi_u$ is either increasing, or decreasing and then increasing, so that $\varphi_u\le\max\{\Phi(0),\varphi_u(t^-(u))\}=\varphi_u(t^-(u))$ there, since $\Phi(0)=0\le c_W^-\le\varphi_u(t^-(u))$ by (H3'). In $(t^-(u),\infty)$ it is either decreasing, or decreasing and then increasing towards a limit which is negative by (H2'). In both cases $\varphi_u<\varphi_u(t^-(u))$ in $(t^-(u),\infty)$, and since $\varphi_u(t^-(u))\ge0$, $\varphi_u$ has exactly one zero $r_1(u)$ in $[t^-(u),\infty)$, with $\varphi_u>0$ in $[t^-(u),r_1(u))$ and $\varphi_u<0$ in $(r_1(u),\infty)$.
\end{proof}

\begin{lemma}\label{l1}
Assume (H1')--(H3'). Then $r_1:W\to(0,\infty)$ is continuous and the set
\begin{equation*}
	A:=\{u\in W: r_1(u)\le 1\}
\end{equation*}
is path-connected. Moreover, $e\in A$ and $\Phi\le 0$ in $A$.
\end{lemma}

\begin{proof}
Let $u_n\to u$ in $W$. Since $\varphi_u(r_1(u)+1)<0$ by Lemma \ref{shape}, for $n$ large we have $\varphi_{u_n}(r_1(u)+1)<0$ and, by Lemma \ref{tcont}, $t^-(u_n)<r_1(u)+1$, so that Lemma \ref{shape} gives $t^-(u_n)\le r_1(u_n)<r_1(u)+1$. Hence every limit point $r$ of $(r_1(u_n))$ satisfies $r\ge t^-(u)$ and $\varphi_u(r)=\lim\Phi(r_1(u_n)u_n)=0$, so that $r=r_1(u)$ by Lemma \ref{shape}. Thus $r_1$ is continuous.

Since $r_1(\sigma u)=\sigma^{-1}r_1(u)$ for $\sigma>0$, we have $A=\{tu: u\in S\cap W,\ t\ge r_1(u)\}$, so that $\Phi\le0$ in $A$ by Lemma \ref{shape}, and $e\in A$, since $t^-(e)<1$ and $\varphi_e(1)=\Phi(e)\le0$. The map $u\mapsto r_1(u)u$ is a homeomorphism from $S\cap W$ onto $B:=\{r_1(u)u: u\in S\cap W\}\subset A$, with inverse $v\mapsto v/\|v\|$, so $B$ is path-connected, and each $tu\in A$ is joined to $r_1(u)u\in B$ by a segment contained in $A$. Hence $A$ is path-connected.
\end{proof}

\begin{proposition}\label{lower}
Assume (H1') and $c_W^-\ge 0$. Then $c_W(e)\geq c^-_W$.
\end{proposition}

\begin{proof}
Since every $\gamma\in\Gamma$ starts at $0$ and $\Phi(0)=0$, we may assume that $c_W^->0$. Let $\rho>0$ be such that $\Phi<c_W^-$ in $B_\rho(0)$, so that $t^-(u)\ge\rho$ for $u\in S\cap W$. The set $O:=\{tu: u\in S\cap W,\ t>t^-(u)\}$ of points lying beyond $\mathcal{N}^-\cap W$ is open, by Lemmas \ref{pathcom} and \ref{tcont}, does not meet $B_\rho(0)$, and contains $e$. Given $\gamma\in\Gamma$, let $s^*:=\inf\{s\in[0,1]: \gamma(s)\in O\}$ and $z:=\gamma(s^*)\in\overline{O}\setminus O$, and write $z=\lim t_nu_n$, with $u_n\in S\cap W$ and $t_n>t^-(u_n)$. Then $\|z\|\ge\rho$ and $u_n\to v:=z/\|z\|\in\overline{W}$. If $v\in W$, then $t^-(u_n)\to t^-(v)\le\|z\|$, and $z\notin O$ forces $\|z\|=t^-(v)$, i.e. $z\in\mathcal{N}^-\cap W$, so that $\Phi(z)\ge c_W^-$. If $v\in\partial W$, then, up to a subsequence, $t^-(u_n)\to\tau\in[\rho,\|z\|]$, so that $\varphi_v(\tau)=\lim\Phi(t^-(u_n)u_n)\ge c_W^->0$, and \eqref{nomax} yields $\Phi(z)=\varphi_v(\|z\|)\ge\varphi_v(\tau)\ge c_W^-$. In both cases $\max_{[0,1]}\Phi\circ\gamma\ge c_W^-$.
\end{proof}

The previous proof does not show that $\gamma$ intersects $\mathcal{N}^-\cap W$, since $z$ may lie in the cone over $\partial W$. This is where (H4') enters.

\begin{proposition}\label{pp2}
Assume (H1'), (H4') and {either $c_W^->0$, or $c_W^-=0$ and (H5') holds}. Then every $\gamma\in\Gamma$ intersects $\mathcal{N}^-\cap W$, and $\mathcal{N}^-\cap W$ is closed in $X$.
\end{proposition}

\begin{proof}
{Let $\rho>0$ be such that $\|v\|\ge\rho$ for every $v\in\mathcal{N}^-\cap W$: if $c_W^->0$ this follows from the continuity of $\Phi$ at $0$, and if $c_W^-=0$ it is (H5').} Let $\gamma\in\Gamma$. Since $W$ is open and $0\notin W\ni e$, we have $s^*:=\max\{s\in[0,1]:\gamma(s)\notin W\}<1$ and $\gamma(s^*)\in\partial W$, and $f(s):=t^-(\gamma(s))$ is continuous in $(s^*,1]$, by Lemma \ref{tcont}. Moreover, $f(s)\to\infty$ as $s\to s^*$: if $\gamma(s^*)=0$, since $f(s)\|\gamma(s)\|\ge\rho$, and otherwise by (H4'), since $f(s)=t^-(\gamma(s)/\|\gamma(s)\|)/\|\gamma(s)\|$ and $\gamma(s)/\|\gamma(s)\|\to\gamma(s^*)/\|\gamma(s^*)\|\in\partial W$. As $f(1)=t^-(e)<1$, we get $f(s_1)=1$, i.e. $\gamma(s_1)\in\mathcal{N}^-\cap W$, for some $s_1\in(s^*,1)$. Finally, if $\mathcal{N}^-\cap W\ni t^-(u_n)u_n\to v$, with $u_n\in S\cap W$, then $\|v\|\ge\rho$, $t^-(u_n)\to\|v\|$ and $u_n\to u:=v/\|v\|\in\overline{W}$. By (H4'), $u\notin\partial W$, so $u\in W$ and, by Lemma \ref{tcont}, $t^-(u)=\|v\|$, i.e. $v\in\mathcal{N}^-\cap W$.
\end{proof}

In order to deal with (H1''), we note that it implies (H4') when $c_W^->0$. {At zero level the same conclusion holds under (H5').}

\begin{lemma}\label{infinity} Suppose (H1'') and {either $c_W^->0$, or $c_W^-=0$ and (H5') holds}. Then (H4') holds.
\end{lemma}
\begin{proof}
Let $(u_n)\subset S\cap W$ be such that $u_n\to u\in\partial W$, and suppose that, along a subsequence, $t^-(u_n)\to t\in[0,\infty)$. Then $\Phi(tu)=\lim\Phi(t^-(u_n)u_n)\ge c_W^-\ge0$, {and $t>0$: if $c_W^->0$ because $\Phi(0)=0$, and if $c_W^-=0$ by (H5').} By (H1''), each $\varphi_{u_n}$ is nonincreasing in $[t^-(u_n),\infty)$, so $\varphi_u$ is nonincreasing in $[t,\infty)$. By \eqref{nomax}, $\varphi_u$ is then constant in $[t,\infty)$, so every point of $(t,\infty)$ is a local maximizer of $\varphi_u$, which is impossible since $u\in\partial W\setminus\{0\}$.
\end{proof}

\subsection{Proofs of Theorems \ref{thmmain} and \ref{thmMPG}}

\begin{proof}[Proof of Theorem \ref{thmmain}]
By Lemma \ref{shape}, $\max_{t>0}\Phi(tu)=\Phi(t^-(u)u)$ for $u\in W$, which gives the second equality. Let $v\in\mathcal{N}^-\cap W$. The segments from $0$ to $v$ and from $v$ to $r_1(v)v$, followed by a path from $r_1(v)v$ to $e$ within $A$ (Lemma \ref{l1}), form a path in $\Gamma$ along which $\Phi\le\Phi(v)$, by Lemmas \ref{shape} and \ref{l1}, since $t^-(v)=1$ and $\Phi(v)\ge c_W^-\ge0$. Thus $c_W(e)\le\Phi(v)$, and taking the infimum over $v$ gives $c_W(e)\le c_W^-$. The reverse inequality is Proposition \ref{lower}.
\end{proof}

\begin{proof}[Proof of Theorem \ref{thmMPG}]
If (H1'') holds, then (H4') holds by Lemma \ref{infinity}. In both cases, Proposition \ref{pp2} shows that $F:=\mathcal{N}^-\cap W\subset W$ is closed in $X$ and meets every $\gamma\in\Gamma$. {Moreover, $\inf_F\Phi=c_W^-\ge0$ and $e\notin F$, since $t^-(e)<1$, so $\Phi$ has the mountain pass geometry with respect to $W$, also when $c_W^-=0$.} Finally, if (H2') holds, then (H1')--(H3') hold, and Theorem \ref{thmmain} yields $\inf_F\Phi=c_W^-=c_W(e)$.
\end{proof}

\begin{remark}
Assume the hypotheses of Theorem \ref{thmmain} and (H4'), and let $\tilde{e}\in W$ be such that $t^-(\tilde{e})>1$ and $\Phi(\tilde{e})\le 0$. Then $c_W(\tilde{e},e)=c_W^-$, where
\begin{equation*}
	c_W(\tilde{e},e):=\inf_{\gamma \in \Gamma(\tilde{e},e)}\max_{t\in [0,1]}\Phi(\gamma(t)), \quad \mbox{and} \quad \Gamma(\tilde{e},e)=\{\gamma\in C([0,1],X): \gamma(0)=\tilde{e},\gamma(1)=e\}.
\end{equation*}
Indeed, as in the proof of Lemma \ref{shape}, $\varphi_{\tilde{e}}\le\max\{0,\Phi(\tilde{e})\}=0$ in $(0,1]$, since $1<t^-(\tilde{e})$, so concatenating the segment from $\tilde{e}$ to $0$ with the paths built in the proof of Theorem \ref{thmmain} gives $c_W(\tilde{e},e)\le c_W^-$. Conversely, a path $\gamma\in\Gamma(\tilde{e},e)$ contained in $W$ meets $\mathcal{N}^-\cap W$, since $t^-(\gamma(0))>1>t^-(\gamma(1))$; otherwise, the proof of Proposition \ref{pp2}, applied after the last time at which $\gamma$ is not in $W$, shows that $\max\Phi\circ\gamma\ge c_W^-$.
\end{remark}

\subsection{Critical points and multiplicity}\label{subcrit}
\begin{corollary}\label{existence}
Let $X$ be a Banach space. Assume the hypotheses of Theorem \ref{thmMPG}, together with (H2'), and suppose that $\Phi\in C^1(X)$ satisfies the Palais--Smale condition at the level $c_W^-$. Then there exists $w\in\mathcal{N}^-\cap W$ such that $\Phi(w)=c_W^-$ and $\Phi'(w)=0$. {In particular, this conclusion remains valid when $c_W^-=0$.}
\end{corollary}

\begin{proof}[Proof of Corollary \ref{existence}]
By Theorem \ref{thmMPG}, $F:=\mathcal{N}^-\cap W$ is closed, meets every path in $\Gamma$, and $\inf_F\Phi=c_W(e)$. Since the components of the open set $X\setminus F$ are path-connected, $F$ separates $0$ and $e$, and the mountain pass principle of Ghoussoub and Preiss \cite{GhPr} yields $(u_n)\subset X$ such that $\mathrm{dist}(u_n,F)\to0$, $\Phi(u_n)\to c_W^-$ and $\Phi'(u_n)\to 0$. By the Palais--Smale condition, a subsequence converges to some $w\in F$ with $\Phi(w)=c_W^-$ and $\Phi'(w)=0$.
\end{proof}

\begin{corollary}\label{multiplicity}
Let $X$ be a Banach space, $\Phi\in C^1(X)$, and let $W_1,\dots,W_k$ be distinct path-connected components of $Y$ such that, for each $j$, the assumptions of Corollary \ref{existence} hold with $W=W_j$. Then $\Phi$ has at least $k$ distinct critical points $w_1,\dots,w_k$, with $w_j\in\mathcal{N}^-\cap W_j$ and $\Phi(w_j)=c^-_{W_j}=c_{W_j}(e_j)$ for every $e_j\in W_j$ such that $\Phi(e_j)\le 0$ and $t^-(e_j)<1$. {The levels $c^-_{W_j}$ may also vanish.}
\end{corollary}

\begin{proof}
Apply Corollary \ref{existence} and Theorem \ref{thmmain} in each $W_j$. The critical points are distinct, since $W_i\cap W_j=\emptyset$ for $i\neq j$.
\end{proof}

\begin{remark}\label{multrem}
Corollary \ref{multiplicity} distinguishes critical points by the component of $Y$ in which they lie, and not by their energy, so that the levels $c^-_{W_j}$ may coincide. The idea of obtaining critical points from different cones is not new: it appears in \cite{FRS}, where a critical point is obtained by minimization in each of two cones, under a compactness condition involving weak limits and the weak boundary of the cones, and in problems where nonnegative and nonpositive solutions are sought separately, by restricting the functional to the cones of nonnegative and nonpositive functions. The contribution of Corollary \ref{multiplicity} is that each of these critical points is a mountain pass point relative to its component, at the level $c^-_{W_j}$, and that only the Palais--Smale condition is required. In Section \ref{secapp} we illustrate it with a planar example having four critical points at the same level (Theorem \ref{thmsat}), with finite-rank problems (Subsection \ref{subfr}), and with a local problem related to the Fu\v{c}\'{\i}k spectrum (Subsection \ref{subfucik}), where the two levels may be equal or different (Proposition \ref{propfucik}).

\end{remark}

Corollary \ref{existence} provides a critical point at the level $c_W^-$ under a compactness condition. Without it, $c_W^-$ need not be achieved. When it is, the minimizer is a critical point of $\Phi$, as the next result shows, and by Theorem \ref{thmmain} it lies at the mountain pass level $c_W(e)$. The first part of the next result follows the argument of \cite[Theorem 2.1]{FRS} (see also \cite{LWW}), which requires neither the differentiability of $t^-$ nor any regularity of the Nehari set.

\begin{proposition}\label{natural}
Let $\Phi\in C(X)\cap C^1(X\setminus\{0\})$ and assume (H1'), (H2') and (H3'). If $w\in\mathcal{N}^-\cap W$ is such that $\Phi(w)=c_W^-$, then $\Phi'(w)=0$. If, moreover, $\Phi$ is of class $C^2$ in a neighborhood of $w$ and $\varphi_w''(1)<0$, then $w$ has Morse index one, i.e. the maximal dimension of a subspace of $X$ on which $\Phi''(w)$ is negative definite is one.
\end{proposition}

\begin{proof}[Proof of Proposition \ref{natural}]
Since $t^-(w)=1$, Lemma \ref{shape} shows that $1$ is a global maximizer of $\varphi_w$. Suppose that $\Phi'(w)v<0$ for some $v\in X$. For $\sigma>0$ small, $w+\sigma v\in W$ and, by Lemma \ref{tcont}, $t_\sigma:=t^-(w+\sigma v)\to1$ as $\sigma\to0$, so that, by the continuity of $\Phi'$ near $w$,
\begin{equation*}
	\Phi(t_\sigma(w+\sigma v))-\Phi(w)\le\Phi(t_\sigma(w+\sigma v))-\Phi(t_\sigma w)=t_\sigma\int_0^\sigma\Phi'(t_\sigma(w+\tau v))v\,d\tau<0
\end{equation*}
for $\sigma$ small. Since $t_\sigma(w+\sigma v)\in\mathcal{N}^-\cap W$, this contradicts $\Phi(w)=c_W^-$. Hence $\Phi'(w)v\ge0$ for every $v\in X$, i.e. $\Phi'(w)=0$.

Assume now that $\Phi$ is of class $C^2$ near $w$ and $\varphi_w''(1)<0$, and let $G(u):=\Phi'(u)u$. Since $\varphi_u'(t)=G(tu)/t$, we have $\varphi_u''(1)=G'(u)u-G(u)$, so that $G'(w)w=\varphi_w''(1)<0$. Given $h\in\ker G'(w)$, the implicit function theorem yields a $C^1$ function $s$ with $s(0)=s'(0)=0$ and $G(u_t)=0$ for $|t|$ small, where $u_t:=w+th+s(t)w$. Then $u_t\in W$, $\varphi_{u_t}'(1)=0>\varphi_{u_t}''(1)=G'(u_t)u_t$, so that $t^-(u_t)=1$ and $\Phi(u_t)\ge c_W^-=\Phi(w)$. Since $\Phi'(w)=0$ and $\Phi(u_t)=\Phi(w)+\frac{t^2}{2}\Phi''(w)[h,h]+o(t^2)$, we obtain $\Phi''(w)[h,h]\ge0$. Hence every subspace on which $\Phi''(w)$ is negative definite meets the hyperplane $\ker G'(w)$ only at $0$, so its dimension is at most one, and $\Phi''(w)[w,w]=\varphi_w''(1)<0$.
\end{proof}

Since $W\subset X\setminus\{0\}$, the differentiability of $\Phi$ at the origin plays no role in Proposition \ref{natural}. This matters for the affine $p$-Laplacian (Theorem \ref{thmaff}), whose energy is only known to be of class $C^1$ away from the origin.

The first assertion of Proposition \ref{natural} only uses the fiber structure and $\Phi\in C^1(X\setminus\{0\})$; in particular, it applies to the concave-convex functional of Theorem \ref{thm5}, which is not of class $C^2$. The Morse index statement applies to the critical points of Theorems \ref{thm4} and \ref{thm6}, whose functionals are of class $C^2$, and where $\varphi_w''(1)$ equals $-(p-2)H_\lambda(w)<0$ and $-H_\lambda(w)<0$, respectively. As for Nehari ground states, index one is the expected index of a critical point of mountain pass type.

\subsection{The classical setting}\label{subclass}
We now consider the situation where all the fiber maps have a local maximizer, which includes {the setting where every fiber map has a unique critical point that is a global maximizer, as in \eqref{equali}}. Let us assume that:
\begin{enumerate}
	\item[(H1)] For any $u\in X\setminus\{0\}$ the map $\varphi_u$ has exactly one local maximizer $t^-(u)$.
	\item[(H2)] $\displaystyle\lim_{t\to\infty}\varphi_u(t)<0$ for every $u\in X\setminus\{0\}$.
\end{enumerate}
Then $Y=X\setminus\{0\}$ and $\mathcal{N}^-=\{t^-(u)u: u\in X\setminus\{0\}\}$, and we set $c^-=\inf_{\mathcal{N}^-}\Phi$. We also assume that:
\begin{enumerate}
	\item[(H3)] $c^-\ge 0$.
\end{enumerate}
Condition (H1) holds, for instance, if every fiber map has a unique critical point, which is a global maximizer, {as in the classical setting discussed in the introduction}, or a local minimizer followed by a local maximizer, as in concave-convex problems with positive weights and small values of the parameter.

\begin{corollary}\label{cor1}
Assume (H1)--(H3) and $\dim X\ge 2$. Then
\begin{equation*}
	c(e)=c^-=\inf_{u\in X\setminus\{0\}}\max_{t>0}\Phi(tu)
\end{equation*}
for every $e\in X$ such that $\Phi(e)\le 0$ and $t^-(e)<1$. If, moreover, $c^->0$, then $\Phi$ has the mountain pass geometry with respect to $X\setminus\{0\}$, with $F=\mathcal{N}^-$.
\end{corollary}

\begin{proof}[Proof of Corollary \ref{cor1}]
By (H1), $W=Y=X\setminus\{0\}$, which is path-connected since $\dim X\ge2$, and (H1)--(H3) imply (H1')--(H3'). Since $\partial W=\{0\}$, (H4') holds trivially, and Theorems \ref{thmmain} and \ref{thmMPG} apply.
\end{proof}

\begin{remark}\label{dimone}
The assumption $\dim X\ge 2$ cannot be dropped. If $\dim X=1$, then $X\setminus\{0\}$ consists of two half-lines $W_\pm$, and Theorem \ref{thmmain} only yields $c(e)=c^-_{W}$, where $W$ is the half-line containing $e$. This value may be strictly larger than $c^-=\min\{c^-_{W_+},c^-_{W_-}\}$.
\end{remark}

\begin{remark}
Corollary \ref{cor1} applies, for instance, to 
\[\Phi(u)=I_1(u)-\lambda I_2(u)-I_3(u), \quad u \in X,\]
where $I_i$ are continuous, positive in $X\setminus\{0\}$ and $\alpha_i$-homogeneous, i.e. $I_i(tu)=t^{\alpha_i}I_i(u)$ for any $t>0$ and $u \in X$, with $0<\alpha_2<\alpha_1<\alpha_3$. In addition, we assume that \[\mu_0:=\inf_{u \in X \setminus \{0\}} \frac{I_1(u)^{\frac{\alpha_3-\alpha_2}{\alpha_3-\alpha_1}}}{I_2(u)I_3(u)^{\frac{\alpha_1-\alpha_2}{\alpha_3-\alpha_1}}}>0.\]
{
Then one can show that (H1)--(H3) hold whenever $0<\lambda<C\mu_0$, for some $C=C(\alpha_1,\alpha_2,\alpha_3)$. In particular, the fiber maxima are positive and $c^-\ge 0$; the assumption $\mu_0>0$ alone does not imply $c^->0$, which would require a uniform positive lower bound for the fiber maxima.

}
 
The condition $\mu_0>0$ is satisfied if, for instance, there exists $C>0$ such that $I_2(u)\leq C I_1(u)^{\alpha_2/\alpha_1}$ and $I_3(u)\leq C I_1(u)^{\alpha_3/\alpha_1}$ for any $u \in  X \setminus \{0\}$.
\end{remark}

\begin{remark}\label{transl}
The same principle holds when dilations are replaced by translations. Let $\Omega\subset\mathbb{R}^N$ be a bounded connected Lipschitz domain, $p>1$, $1<q<p$ and $a\in L^\infty(\Omega)$ with $a\ge a_*>0$, and consider the Neumann problem $-\Delta_p u=a(x)|u|^{q-2}u$ in $\Omega$, $|\nabla u|^{p-2}\partial_\nu u=0$ on $\partial\Omega$, whose energy is
\begin{equation*}
	\Phi(u)=\frac{1}{p}\int_\Omega|\nabla u|^p-\frac{1}{q}\int_\Omega a|u|^q, \quad u\in W^{1,p}(\Omega).
\end{equation*}
Since the gradient term is invariant under the addition of constants, each fiber $c\mapsto \Phi(u+c)$ is strictly concave and tends to $-\infty$ as $|c|\to\infty$. Hence the set
\begin{equation*}
	\mathcal{N}_f:=\Big\{u\in W^{1,p}(\Omega): \int_\Omega a|u|^{q-2}u=0\Big\}=\{u\in W^{1,p}(\Omega): \Phi'(u)1=0\}
\end{equation*}
meets every line $\{u+c: c\in\mathbb{R}\}$ exactly once, at the maximum of the fiber. For $p=2$ and $a\equiv 1$ this constant-shift reduction is due to Parini and Weth \cite{PaWe}. Since $\int_\Omega a|u|^{q-2}u$ is negative at the constant function $-R$ and positive at $R$, every path joining $-R$ to $R$ in $W^{1,p}(\Omega)$ intersects $\mathcal{N}_f$. Arguing as in the proof of Theorem \ref{thmmain}, with paths made of the segment from $-R$ to $u-R$, the fiber $\{u+c: |c|\le R\}$ through $u\in\mathcal{N}_f$, and the segment from $u+R$ to $R$, and using that minimizing sequences of $\Phi$ on $\mathcal{N}_f$ are bounded, one obtains
\begin{equation*}
	\inf_{\mathcal{N}_f}\Phi=\inf_{\gamma}\max_{t\in[0,1]}\Phi(\gamma(t)) \quad \mbox{for $R$ large},
\end{equation*}
where $\gamma$ ranges over the continuous paths from $-R$ to $R$ in $W^{1,p}(\Omega)$. This is the translation analogue of Corollary \ref{cor1}: every fiber has a maximum, so $\mathcal{N}_f$ is a graph over the space of functions with zero mean, and no cone is needed. Note, however, that $0\in\mathcal{N}_f$ and that $\inf_{\mathcal{N}_f}\Phi<0$: the roles of $0$ and $e$ are played by the constants $-R$ and $R$, at which $\Phi$ lies far below $\inf_{\mathcal{N}_f}\Phi$.
\end{remark}

\section{Examples and applications}\label{secapp}
Each subsection below starts with the setting and the statement of the corresponding result, followed by its proof.

\subsection{Finite dimensional examples}\label{subfin}
We start with a finite dimensional example, which shows that it is possible to find different cones for which we can apply Theorem \ref{thmmain}. Let $\Phi:\mathbb{R}^2\to \mathbb{R}$ be defined by 
\begin{equation*}
	\Phi(x,y)=\frac{1}{2}g\left(\frac{(x,y)}{\|(x,y)\|}\right)(x^2+y^2)-\frac{xy}{6}(x^2+y^2)^2, (x,y)\in \mathbb{R}^2\setminus \{0\},\ \ \ \Phi(0,0)=0,
\end{equation*}
where $g:S\to \mathbb{R}$ is a $C^1$ function defined over the sphere $S=\{(x,y)\in \mathbb{R}^2: \|(x,y)\|=1\}$. Clearly $\Phi$ is $C^1$ and $\Phi'(0)=0$. Let $W_1,W_2,W_3,W_4$ be the first, second, third and fourth open quadrants of $\mathbb{R}^2$, respectively.

\begin{theorem}\label{thm3}
Assume that $g>0$ on $S$. Then $Y=W_1\cup W_3$, and (H1')--(H4') and (H1'') hold for $W=W_1$ and for $W=W_3$, with $c_{W_1}^-,c_{W_3}^->0$. Moreover, for $i=1,3$, we have $c_{W_i}(e_i)=c_{W_i}^-$ for every $e_i\in W_i$ such that $\Phi(e_i)<0$, and there exists $w_i\in\mathcal{N}^-\cap W_i$ such that $\Phi(w_i)=c_{W_i}^-$ and $\Phi'(w_i)=0$.
\end{theorem}
\begin{remark}\label{remg} For a suitable choice of $g$ in Theorem \ref{thm3}, we have $c_{W_1}^-<c_{W_3}^-$. For instance, let $g(x,y)=2-x-y$, so that $g>0$ on $S$. By \eqref{finiex}, $c_{W_i}^-=\inf_{S\cap W_i}\frac{g^{3/2}}{3(xy)^{1/2}}$, and this infimum is achieved, since the quotient tends to $\infty$ at $\partial W_i$. The map $u\mapsto -u$ takes $S\cap W_3$ onto $S\cap W_1$, preserves $xy$ and satisfies $g(-u)<2<g(u)$ for $u\in W_3$. Hence $c_{W_1}^-<c_{W_3}^-$. Moreover, if
	\begin{equation*}
		c=\inf_{\gamma \in \Gamma_0}\sup _{t\in [0,1]}\Phi(\gamma(t)),
	\end{equation*}
	where  $\Gamma_0=\{\gamma\in C([0,1],\mathbb{R}^2): \gamma(0)=0,\ \Phi(\gamma(1))<0\}$, then $c=c_{W_1}^-<c_{W_3}^-$ for this choice of $g$. Indeed, $\Phi>0$ in $\overline{W_2\cup W_4}\setminus\{0\}$, so that $\gamma(1)\in W_1\cup W_3$ for every $\gamma\in\Gamma_0$, and $c_{W_i}(e)=c_{W_i}^-$ for every $e\in W_i$ with $\Phi(e)<0$, by Theorem \ref{thmmain}. 
\end{remark}
\begin{remark}\label{remh}
More generally, let $g>0$ on $S$, let $h$ be a homogeneous polynomial of degree $d\ge1$ which is positive somewhere, and let
\begin{equation*}
	\Phi(u)=\frac{1}{2}g\left(\frac{u}{\|u\|}\right)\|u\|^2-h(u)\|u\|^4, \quad u\in\mathbb{R}^2\setminus\{0\}, \qquad \Phi(0)=0.
\end{equation*}
For $u\in S$ we have $\varphi_u(t)=\frac{1}{2}g(u)t^2-h(u)t^{d+4}$, so that $Y=\{h>0\}$, and for $u\in S\cap Y$,
\begin{equation*}
	t^-(u)=\left(\frac{g(u)}{(d+4)h(u)}\right)^{\frac{1}{d+2}} \quad \mbox{and} \quad \varphi_u(t^-(u))=\frac{d+2}{2(d+4)}\,g(u)\,t^-(u)^2.
\end{equation*}
Arguing as in the proof of Theorem \ref{thm3}, one sees that its conclusions hold in every path-connected component of $Y$. For instance, if $h(x,y)=xy(x^2-y^2)$, which equals $\frac{1}{4}r^4\sin 4\theta$ in polar coordinates, then $Y$ consists of four open sectors of angle $\frac{\pi}{4}$, each of which contains a critical point of $\Phi$ at its own mountain pass level.
\end{remark}

By Corollary \ref{cor1}, under (H1)--(H3) with $c^->0$ the set $\mathcal{N}^-$ separates $0$ from $e$, so that $\Phi$ satisfies the mountain pass geometry in the sense of Ghoussoub and Preiss \cite[Theorem 1]{GhPr} with $F=\mathcal{N}^-$ (see also Ambrosetti et al. \cite{ABC}). We now give an example where $\Phi$ does not satisfy the classical mountain pass geometry; in fact, $\Phi$ has no local minimizers.

\begin{theorem}\label{thm3.1}
Assume that $g>0$ on $\overline{W_1\cup W_2\cup W_4}\cap S$ and that $g(x_0,y_0)<0$ for some $(x_0,y_0)\in W_3\cap S$. Then
\begin{equation*}
	Y=W_1\cup\left\{u\in W_3: g\left(\frac{u}{\|u\|}\right)>0\right\},
\end{equation*}
$W_1$ is a path-connected component of $Y$, and (H1')--(H4') and (H1'') hold for $W=W_1$, with $c_{W_1}^->0$. Moreover, $c_{W_1}(e_1)=c_{W_1}^-$ for every $e_1\in W_1$ such that $\Phi(e_1)<0$, and $\Phi$ has no local minimizers.
\end{theorem}

\begin{figure}[ht]
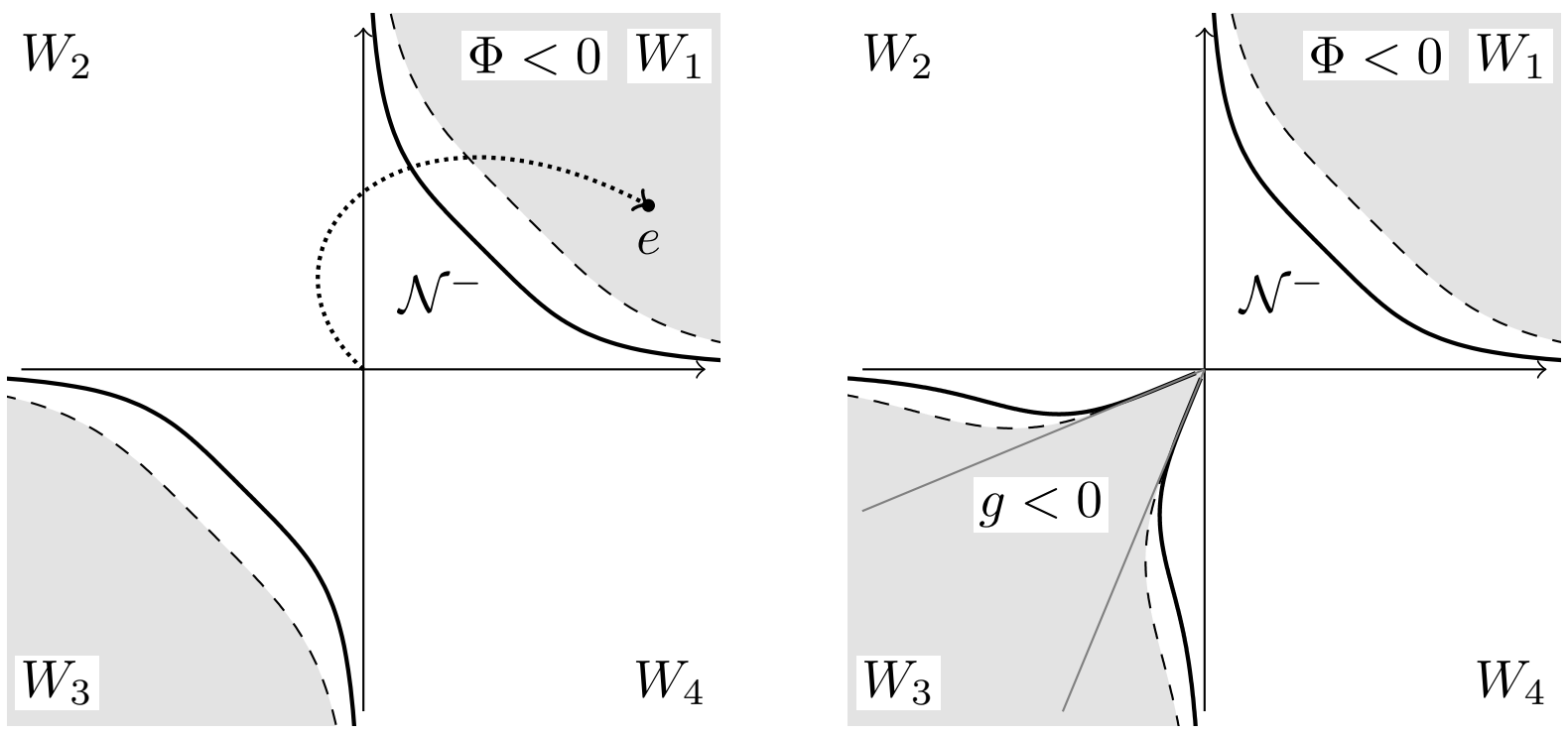

\centering

\caption{The set $\mathcal{N}^-$ (solid), the curve $\{\Phi=0\}$ (dashed) and the region $\{\Phi<0\}$ (shaded) for the functional of Theorems \ref{thm3} and \ref{thm3.1}. Left: $g\equiv1$; the dotted path joins $0$ to $e$, leaves $W_1$ and still has to cross $\mathcal{N}^-\cap W_1$. Right: $g(x,y)=1-8\min\{x,0\}^2\min\{y,0\}^2$, which is negative in a sector of $W_3$. There $\Phi<0$ along every ray, so $0$ is not a local minimizer, while $\mathcal{N}^-\cap W_1$ still separates $0$ from every $e\in W_1$ with $\Phi(e)<0$. The cones of $Y$ inside $W_3$ have $c_W^-=0$, since $\mathcal{N}^-$ reaches the origin there.}
\label{fig1}
\end{figure}

\begin{remark}\label{MPGC}
By Theorem \ref{thmMPG}, the function $\Phi$ in Theorem \ref{thm3.1} has the mountain pass geometry with respect to $W_1$, with $F=\mathcal{N}^-\cap W_1$. On the other hand, it does not satisfy the classical mountain pass geometry: if $\Phi\ge\alpha>0$ on some sphere $\partial B_\rho$, then, since $\Phi(0)=0$, the minimum of $\Phi$ over the compact ball $\overline{B_\rho}\subset\mathbb{R}^2$ would be achieved at an interior point, which would be a local minimizer of $\Phi$. See Figure \ref{fig1}. A PDE example with the same features is given in Theorem \ref{thm6}.

\end{remark}

The cones $W_1$ and $W_3$ in Theorem \ref{thm3} are full quadrants because the nonlinear term is homogeneous. For a nonlinear term that is only asymptotically quadratic, the cones are determined by the behavior of the fiber maps at infinity. Consider the smooth function
\begin{equation}\label{satur}
	\Phi(x,y)=\frac{1}{2}(x^2+y^2)-\frac{8x^2y^2}{1+x^2+y^2}, \quad (x,y)\in\mathbb{R}^2,
\end{equation}
which satisfies $\Phi(0)=0$ and $\Phi'(0)=0$, and let
\begin{equation*}
	\mathcal{W}=\left\{(x,y)\in\mathbb{R}^2: (2-\sqrt3)|y|<|x|<(2+\sqrt3)|y|\right\}=\left\{(x,y)\in\mathbb{R}^2: (x^2+y^2)^2<16x^2y^2\right\}.
\end{equation*}
Then $\mathcal{W}$ has four connected components $\mathcal{W}_1,\dots,\mathcal{W}_4$, one in each open quadrant, each of them bisected by one of the diagonals $y=\pm x$.

\begin{theorem}\label{thmsat}
Let $\Phi$ be given by \eqref{satur}. Then $Y=\mathcal{W}$, and (H1''), (H2') and (H3') hold in each $\mathcal{W}_i$, with
\begin{equation*}
	c^-_{\mathcal{W}_i}=\frac{7}{2}-2\sqrt3>0, \quad i=1,\dots,4.
\end{equation*}
In particular, $c_{\mathcal{W}_i}(e_i)=c^-_{\mathcal{W}_i}$ for every $e_i\in\mathcal{W}_i$ such that $\Phi(e_i)\le 0$ and $t^-(e_i)<1$, for instance $e_i=(\pm2,\pm2)$, and $\Phi$ has the mountain pass geometry with respect to each $\mathcal{W}_i$. Moreover, $c^-_{\mathcal{W}_i}$ is achieved at the critical points $(\pm s,\pm s)$ of $\Phi$, where $s=\sqrt{(2\sqrt3-3)/6}$.
\end{theorem}

\begin{remark}
In contrast with Remark \ref{remg}, where the levels $c^-_{W_1}$ and $c^-_{W_3}$ can be made different, here $\Phi$ is invariant under $x\mapsto -x$, $y\mapsto -y$ and $(x,y)\mapsto(y,x)$, so the four levels coincide. Since $\Phi\ge 0$ outside $\mathcal{W}$ (see the proof below), every path from $0$ to a point where $\Phi<0$ ends in some $\mathcal{W}_i$. Hence the four points $(\pm s,\pm s)$ are critical points at the level $\inf_{\gamma\in\Gamma_0}\max_{t\in[0,1]}\Phi(\gamma(t))$, with $\Gamma_0$ as in Remark \ref{remg}, and they are distinguished only by the component of $Y$ in which they lie (see Corollary \ref{multiplicity}). The cones $\mathcal{W}_i$ are proper subsectors of the quadrants because the nonlinear term is asymptotically quadratic: along every ray it competes with $\frac{1}{2}(x^2+y^2)$, and it prevails only inside $\mathcal{W}$.
\end{remark}

Let us now prove Theorems \ref{thm3}, \ref{thm3.1} and \ref{thmsat}. For the functional of Theorems \ref{thm3} and \ref{thm3.1} and $u=(x,y)\in S$, we have
\begin{equation*}
	\varphi_u(t)=\frac{g(u)}{2}t^2-\frac{xy}{6}t^6 \quad \mbox{and} \quad \varphi_u'(t)=t\left(g(u)-xy\,t^4\right), \qquad t>0.
\end{equation*}
Under the assumptions of both theorems, $g>0$ on the coordinate axes, so $g(u)$ and $xy$ do not vanish simultaneously. Hence $\varphi_u$ has a critical point if, and only if, $g(u)$ and $xy$ are both positive or both negative. In the first case the critical point is unique and is a global maximizer, given by
\begin{equation}\label{finiex}
	t^-(u)=\left[\frac{g(u)}{xy}\right]^{\frac{1}{4}} \quad \mbox{and} \quad \varphi_u(t^-(u))=\frac{1}{3}\frac{g(u)^{3/2}}{(xy)^{1/2}},
\end{equation}
and $\varphi_u(t)\to-\infty$ as $t\to\infty$; in the second case it is a global minimizer. Thus $Y=\{(x,y)\neq0: xy>0,\ g((x,y)/\|(x,y)\|)>0\}$, and (H1'), (H1'') and (H2') hold in every path-connected component of $Y$. Moreover, if $e\in Y$ and $\Phi(e)<0$, then $t^-(e)<1$, since $\varphi_e$ is positive in $(0,t^-(e)]$.

\begin{proof}[Proof of Theorem \ref{thm3}]
Since $g>0$ on $S$, $Y=W_1\cup W_3$. Let $i\in\{1,3\}$. Since $0<xy\le\frac{1}{2}$ on $S\cap W_i$, \eqref{finiex} gives $\varphi_u(t^-(u))\ge\frac{\sqrt2}{3}\left(\min_S g\right)^{3/2}$ for $u\in S\cap W_i$, and $t^-(u),\varphi_u(t^-(u))\to\infty$ as $u\in S\cap W_i$ approaches $\partial W_i$. Thus (H3') and (H4') hold, with $c_{W_i}^->0$, Theorem \ref{thmmain} gives $c_{W_i}(e_i)=c_{W_i}^-$, and $c_{W_i}^-=\inf_{S\cap W_i}\varphi_u(t^-(u))$ is achieved at some $w_i=t^-(u_i)u_i$, which is a critical point by Proposition \ref{natural}.
\end{proof}

\begin{proof}[Proof of Theorem \ref{thm3.1}]
The description of $Y$ follows from the discussion above, and $W_1$ is a component of $Y$ since $\overline{W_1}\cap\overline{W_3}=\{0\}$. The assertions on $W_1$ follow as in the proof of Theorem \ref{thm3}, with $\min_{\overline{W_1}\cap S}g$ instead of $\min_S g$. Since $g(u_0)<0<x_0y_0$, we have $\Phi<0$ on the ray of $u_0=(x_0,y_0)$, so $0$ is not a local minimizer. If $u\neq0$ is a critical point of $\Phi$, then $\varphi_u'(1)=0$, and since $g>0$ on $\overline{W_2\cup W_4}\cap S$, $g(u/\|u\|)$ and $xy$ are positive, so that $1$ is the strict global maximizer of $\varphi_u$, and $u$ is not a local minimizer either. Finally, if $W$ is a component of $\{u\in W_3: g(u/\|u\|)>0\}$, then by \eqref{finiex}, $t^-(u)\to0$ and $\varphi_u(t^-(u))\to0$ as $u\in S\cap W$ approaches $\partial W\cap W_3$, where $g$ vanishes, so that $c_W^-=0$ and (H4') fails; still, Theorem \ref{thmmain} yields $c_W(e)=c_W^-=0$.
\end{proof}

\begin{proof}[Proof of Theorem \ref{thmsat}]
Let $u=(x,y)$ with $x^2+y^2=1$, and set $k=8x^2y^2\in[0,2]$ and $\tau=t^2$. Then
\begin{equation*}
	\varphi_u(t)=\frac{\tau}{2}-\frac{k\tau^2}{1+\tau}, \qquad \varphi_u'(t)=t\left(1-k\,h(\tau)\right), \quad \mbox{where} \quad h(\tau)=\frac{2\tau^2+4\tau}{(1+\tau)^2}.
\end{equation*}
Since $h'(\tau)=4(1+\tau)^{-3}>0$, $h(0)=0$ and $h(\infty)=2$, $\varphi_u$ has a critical point if, and only if, $2k>1$, i.e. $u\in\mathcal{W}$; it is then unique and a global maximizer, $\tau^-:=t^-(u)^2$ satisfies $(1+\tau^-)^2=\frac{2k}{2k-1}$, and $\varphi_u(t)\to-\infty$ as $t\to\infty$. If $2k\le 1$, then $\varphi_u(t)\ge\frac{\tau}{2}-k\tau\ge 0$; in particular $\Phi\ge0$ outside $\mathcal{W}$. Thus $Y=\mathcal{W}$, and (H1'') and (H2') hold in each $\mathcal{W}_i$. Since $\varphi_u(t)$ decreases in $k$, so does $E(k):=\max_{t>0}\varphi_u(t)>0$, and $c^-_{\mathcal{W}_i}$ is achieved on the diagonals, where $k=2$, $\tau^-=\frac{2}{\sqrt3}-1$ and
\begin{equation*}
	c^-_{\mathcal{W}_i}=E(2)=\frac{\tau^-}{2}-\frac{2(\tau^-)^2}{1+\tau^-}=\frac{7}{2}-2\sqrt3>0,
\end{equation*}
i.e. at the points $(\pm s,\pm s)$, $s^2=\tau^-/2=(2\sqrt3-3)/6$, which are critical points by Proposition \ref{natural}. Finally, $\Phi(\pm2,\pm2)=4-\frac{128}{9}<0$, and the remaining assertions follow from Theorems \ref{thmmain} and \ref{thmMPG}.
\end{proof}

\subsection{An indefinite type functional with two homogeneous terms}\label{subind}
From now on, $\Omega$ is a bounded domain of $\mathbb{R}^N$, $2^*=\frac{2N}{N-2}$ if $N\ge3$ and $2^*=\infty$ if $N=1,2$, and $\|u\|=\left(\int_\Omega|\nabla u|^2\right)^{1/2}$ is the norm of $H^1_0(\Omega)$. We denote by $\lambda_1$ the first eigenvalue of $-\Delta$ in $H^1_0(\Omega)$ and by $\phi_1>0$ an associated eigenfunction. Let $2<p<2^*$ and let $b\in L^\infty(\Omega)$, which may change sign, with $b^+\not\equiv0$. The weak solutions of
\begin{equation*}
	-\Delta u=\lambda u+b(x)|u|^{p-2}u \quad \mbox{in } \Omega, \qquad u=0 \quad \mbox{on } \partial\Omega,
\end{equation*}
are the critical points of
\begin{equation*}
	\Phi_\lambda(u)=\frac{1}{2}H_\lambda(u)-\frac{1}{p}F(u), \quad \mbox{where} \quad H_\lambda(u)=\int_\Omega\left(|\nabla u|^2-\lambda u^2\right) \quad \mbox{and} \quad F(u)=\int_\Omega b|u|^p.
\end{equation*}
Let
\begin{equation*}
	\lambda^*=\inf\left\{\frac{\|u\|^2}{\int_\Omega u^2}: u\in H^1_0(\Omega),\ F(u)>0\right\}.
\end{equation*}
Clearly $\lambda^*\ge\lambda_1$. Moreover, $\lambda^*>\lambda_1$ if, and only if, $\int_\Omega b\,\phi_1^p<0$ (Lemma \ref{li3}).

\begin{theorem}\label{thm4}
Let $\lambda<\lambda^*$. Then $Y=\{u\in H^1_0(\Omega): F(u)>0\}$, and:
\begin{enumerate}[label=(\roman*)]
	\item There exists $c_0>0$ such that, for every path-connected component $W$ of $Y$ and every $e\in W$ with $\Phi_\lambda(e)\le0$ and $t^-(e)<1$,
	\begin{equation*}
		c_W(e)=c_W^-=\inf_{u\in W}\max_{t>0}\Phi_\lambda(tu)\ge c_0,
	\end{equation*}
	and $\Phi_\lambda$ has the mountain pass geometry with respect to $W$.
	\item There exist a path-connected component $W_\lambda$ of $Y$ and $w_\lambda\in\mathcal{N}^-\cap W_\lambda$ such that $\Phi_\lambda(w_\lambda)=c^-_{W_\lambda}=\inf_{\mathcal{N}^-}\Phi_\lambda$ and $\Phi_\lambda'(w_\lambda)=0$.
\end{enumerate}
\end{theorem}

The existence of a critical point at the level $\inf_{\mathcal{N}^-}\Phi_\lambda$ is known: see \cite{RS}, and \cite{FRS}, where ground states relative to cones are obtained for problems of this type. The new information in Theorem \ref{thm4} is the minimax characterization of $c_W^-$ on each path-connected component $W$ of $Y$. The next lemma is essentially contained in \cite{RS}; we include short proofs of it and of Theorem \ref{thm4}(ii) for completeness.

\begin{lemma}\label{li3}
\begin{enumerate}[label=(\roman*)]
	\item $\lambda^*=\inf\left\{\frac{\|u\|^2}{\int_\Omega u^2}: u\in H^1_0(\Omega)\setminus\{0\},\ F(u)\ge0\right\}$.
	\item $\lambda^*>\lambda_1$ if, and only if, $\int_\Omega b\,\phi_1^p<0$.
	\item If $\lambda<\lambda^*$, then there exists $C>0$ such that $H_\lambda(u)\ge C\|u\|^2$ for every $u\in H^1_0(\Omega)$ such that $F(u)\ge0$.
\end{enumerate}
\end{lemma}

\begin{proof}
Let $R(u)=\|u\|^2/\int_\Omega u^2$ for $u\neq0$, and recall that $F\in C^1(H^1_0(\Omega))$ is weakly continuous.

(i) It suffices to approximate every $u\neq0$ with $F(u)=0$ by functions $w$ with $F(w)>0$. If $F'(u)v>0$ for some $v$, take $w=u+\varepsilon v$, $\varepsilon\to0^+$. Otherwise $b|u|^{p-2}u=0$ a.e., so that $F(u+\varepsilon v)=\varepsilon^pF(v)$, and we take $v$ with $F(v)>0$, which exists since $b^+\not\equiv0$.

(ii) If $\int_\Omega b\,\phi_1^p\ge0$, then (i) yields $\lambda^*\le R(\phi_1)=\lambda_1$. Conversely, if $\lambda^*=\lambda_1$, let $\|u_n\|=1$, $F(u_n)>0$, $R(u_n)\to\lambda_1$ and $u_n\rightharpoonup u$. Then $\lambda_1\int_\Omega u^2=1\ge\|u\|^2$, so that $u\neq0$ and $R(u)\le\lambda_1$, i.e. $u=s\phi_1$ with $s\neq0$, and $|s|^p\int_\Omega b\,\phi_1^p=F(u)\ge0$.

(iii) Otherwise there are $u_n\rightharpoonup u$ with $\|u_n\|=1$, $F(u_n)\ge0$ and $H_\lambda(u_n)\to\ell\le0$. Then $F(u)\ge0$ and $\lambda\int_\Omega u^2=1-\ell\ge1\ge\|u\|^2$, so that $u\neq0$ and $R(u)\le\lambda<\lambda^*$, which contradicts (i).
\end{proof}

\begin{proof}[Proof of Theorem \ref{thm4}]
(i) For $u\neq0$ we have $\varphi_u(t)=\frac{t^2}{2}H_\lambda(u)-\frac{t^p}{p}F(u)$ and $\varphi_u'(t)=t\left(H_\lambda(u)-t^{p-2}F(u)\right)$. Hence $\varphi_u$ has a local maximizer if, and only if, $H_\lambda(u)>0$ and $F(u)>0$. In this case it is the unique critical point of $\varphi_u$, given by $t^-(u)=\left(H_\lambda(u)/F(u)\right)^{\frac{1}{p-2}}$, and
\begin{equation}\label{energynehari}
	\varphi_u(t^-(u))=\frac{p-2}{2p}\,\frac{H_\lambda(u)^{\frac{p}{p-2}}}{F(u)^{\frac{2}{p-2}}}.
\end{equation}
By Lemma \ref{li3}(iii), $Y=\{F>0\}$, which is open, and (H1'') and (H2') hold in every component of $Y$. Since $F(u)\le C_1\|u\|^p$, Lemma \ref{li3}(iii) and \eqref{energynehari} yield $\varphi_u(t^-(u))\ge\frac{p-2}{2p}C^{\frac{p}{p-2}}C_1^{-\frac{2}{p-2}}=:c_0>0$ for $u\in Y$, so (H3') holds, and Theorems \ref{thmmain} and \ref{thmMPG} yield (i).

(ii) Let $J(u)=\varphi_u(t^-(u))$, so that $c:=\inf_{\mathcal{N}^-}\Phi_\lambda=\inf_YJ$, and let $(u_n)\subset Y$ be such that $\|u_n\|=1$, $J(u_n)\to c$ and $u_n\rightharpoonup u$. Since $H_\lambda(u_n)\ge C$ and $(J(u_n))$ is bounded, \eqref{energynehari} shows that $F(u)=\lim F(u_n)>0$, so that $u\in Y$ and, by weak lower semicontinuity, $c\le J(u)\le\liminf J(u_n)=c$. Thus $w_\lambda:=t^-(u)u$ satisfies $\Phi_\lambda(w_\lambda)=c\le c^-_{W_\lambda}\le\Phi_\lambda(w_\lambda)$, where $W_\lambda$ is the component of $Y$ containing $u$, and $\Phi_\lambda'(w_\lambda)=0$ by Proposition \ref{natural}.
\end{proof}

\subsection{A concave-convex type problem}\label{subcc}
We now consider an example where the fiber maps may have two critical points. Let $1<q<2<p<2^*$ and let $a,b\in L^\infty(\Omega)$, which may change sign, with $b^+\not\equiv0$. The weak solutions of
\begin{equation*}
	-\Delta u=\lambda a(x)|u|^{q-2}u+b(x)|u|^{p-2}u \quad \mbox{in } \Omega, \qquad u=0 \quad \mbox{on } \partial\Omega,
\end{equation*}
are the critical points of
\begin{equation*}
	\Phi_\lambda(u)=\frac{1}{2}H(u)-\frac{\lambda}{q}G(u)-\frac{1}{p}F(u),
\end{equation*}
where $H(u)=\|u\|^2$, $G(u)=\int_\Omega a|u|^q$ and $F(u)=\int_\Omega b|u|^p$.
Note that $\Phi_\lambda\in C^1(H^1_0(\Omega))$, but $\Phi_\lambda$ is not of class $C^2$ in general, since $q<2$. Let $V=\{u\in H_0^1(\Omega): G(u)>0,\ F(u)>0\}$, which we assume to be nonempty, and
\begin{equation*}
\lambda_0(u)=\frac{q(p-2)}{2(p-q)}\left(\frac{p(2-q)}{2(p-q)}\right)^{\frac{2-q}{p-2}}\frac{H(u)^{\frac{p-q}{p-2}}}{G(u)F(u)^{\frac{2-q}{p-2}}},\quad u\in V.
\end{equation*}
Note that $\lambda_0(u)$ is the unique $\lambda>0$ for which $\varphi_u$ has a critical point at the zero energy level (see the proof of Lemma \ref{lbd0}).

By the Sobolev embedding, $G(u)\le CH(u)^{q/2}$ and $F(u)\le CH(u)^{p/2}$, and since $\frac{q}{2}+\frac{p}{2}\cdot\frac{2-q}{p-2}=\frac{p-q}{p-2}$, it follows that $\lambda_0^*:=\inf_{u\in V}\lambda_0(u)>0$.

\begin{theorem}\label{thm5}
Let $\lambda\in(0,\lambda_0^*]$. Then $Y=\{u\in H^1_0(\Omega): F(u)>0\}$, and:
\begin{enumerate}[label=(\roman*)]
	\item For every path-connected component $W$ of $Y$ and every $e\in W$ with $\Phi_\lambda(e)\le0$ and $t^-(e)<1$,
	\begin{equation*}
		c_W(e)=c_W^-=\inf_{u\in W}\max_{t>0}\Phi_\lambda(tu).
	\end{equation*}
	If $\lambda<\lambda_0^*$, then $c_W^->0$ and $\Phi_\lambda$ has the mountain pass geometry with respect to $W$. {If $\lambda=\lambda_0^*$, the same geometry holds with respect to $W$, including whenever $c_W^-=0$.}
	\item There exist a path-connected component $W_\lambda$ of $Y$ and $w_\lambda\in\mathcal{N}^-\cap W_\lambda$ such that $\Phi_\lambda(w_\lambda)=c^-_{W_\lambda}=\inf_{\mathcal{N}^-}\Phi_\lambda$ and $\Phi_\lambda'(w_\lambda)=0$. Moreover, $c^-_{W_\lambda}>0$ if $\lambda<\lambda_0^*$, while $c^-_{W_\lambda}=0$ if $\lambda=\lambda_0^*$. {At $\lambda=\lambda_0^*$, $\Phi_\lambda$ has zero-level mountain pass geometry with respect to $W_\lambda$.}
\end{enumerate}
\end{theorem}

\begin{remark}\label{zerolevel}
For $\lambda=\lambda_0^*$, the mountain pass level $c_{W_\lambda}(e)$ vanishes, i.e. it equals $\Phi_\lambda(0)$, and still it is achieved by the nontrivial critical point $w_{\lambda_0^*}$ (see also \cite{RSS} on zero energy critical points). {This is a zero-level instance of Definition \ref{defMPG}: $\mathcal{N}^-\cap W_{\lambda_0^*}$ is closed in $X$, separates $0$ from every admissible endpoint $e\in W_{\lambda_0^*}$, and its infimum energy is zero. Zero-altitude mountain pass situations were studied by Pucci and Serrin {\cite{PuSe85}} and by Ghoussoub and Preiss {\cite{GhPr}}. Here the equality between the Nehari infimum and the minimax level is obtained from Theorem \ref{thmmain} without any compactness assumption.}
\end{remark}

Let us compare Theorem \ref{thm5} with known results. If $u\in V$, then $\varphi_u$ has a local maximizer if, and only if, $\lambda<\Lambda(u)$, where
\begin{equation*}
	\Lambda(u):=\frac{p-2}{p-q}\left(\frac{2-q}{p-q}\right)^{\frac{2-q}{p-2}}\frac{H(u)^{\frac{p-q}{p-2}}}{G(u)F(u)^{\frac{2-q}{p-2}}}
\end{equation*}
(see Lemma \ref{lbd0}), and $\lambda_0(u)=\frac{q}{2}\left(\frac{p}{2}\right)^{\frac{2-q}{p-2}}\Lambda(u)<\Lambda(u)$, since $\log\frac{q}{2}<-\frac{2-q}{2}$ and $\log\frac{p}{2}<\frac{p-2}{2}$. The map $\Lambda$ is the nonlinear Rayleigh quotient of Il'yasov \cite{I1}. For $0<\lambda<\inf_{V}\Lambda$, every $u\in V$ has a fiber map with exactly two critical points, and \cite[Corollary 3.2]{FRS} (see also \cite{BW,I1,Ta}) provides a local minimizer of $\Phi_\lambda$ at a negative level and a second critical point, associated with the second critical point of the fiber maps. For $\lambda>\lambda_0^*$, the energy of the latter may be negative; then (H3') fails and Theorem \ref{thmmain} does not apply. For $\lambda\le\lambda_0^*$, the existence of $w_\lambda$ is known \cite{RSiS} (see also \cite{RSS} for $\lambda=\lambda_0^*$), and Theorem \ref{thm5} adds the minimax characterization of its level. We include a short proof of Theorem \ref{thm5}(ii) for completeness.

\begin{lemma}\label{lbd0}
Let $\lambda>0$ and $u\in H^1_0(\Omega)\setminus\{0\}$.
\begin{enumerate}[label=(\roman*)]
	\item The map $\varphi_u$ has at most two critical points. If it has two, the first one is a local minimizer and the second one is a local maximizer.
	\item The map $\varphi_u$ has a local maximizer if, and only if, either $F(u)>0$ and $G(u)\le0$, or $u\in V$ and $\lambda<\Lambda(u)$. In this case $\varphi_u(t)\to-\infty$ as $t\to\infty$.
	\item If $\lambda\le\lambda_0^*$ and $F(u)>0$, then $\varphi_u$ has a local maximizer $t^-(u)$, which is a global maximizer of $\varphi_u$, and $\varphi_u(t^-(u))\ge0$, with strict inequality if $\lambda<\lambda_0^*$. Moreover, if $u\in V$ and $\lambda=\lambda_0(u)$, then $\varphi_u(t^-(u))=0$.
\end{enumerate}
\end{lemma}

\begin{proof}
We have $\varphi_u'(t)=t^{q-1}\psi_u(t)$, where $\psi_u(t)=t^{2-q}H(u)-t^{p-q}F(u)-\lambda G(u)$.

(i)--(ii) If $F(u)\le0$, then $\psi_u$ is increasing, so $\varphi_u$ has at most one critical point, a local minimizer. If $F(u)>0$, then $\psi_u$ increases in $(0,s_u]$, decreases in $[s_u,\infty)$ and tends to $-\infty$, where $s_u^{p-2}=\frac{(2-q)H(u)}{(p-q)F(u)}$. Thus $\psi_u$ has at most two zeros, $\varphi_u$ has a local minimizer at the first one and a local maximizer at the second one if there are two, and $\varphi_u$ has a local maximizer if, and only if, $\psi_u(s_u)>0$. Since a direct computation gives $\psi_u(s_u)+\lambda G(u)=\max_{t>0}\left(t^{2-q}H(u)-t^{p-q}F(u)\right)>0$, with the maximum equal to $\Lambda(u)G(u)$ if $G(u)>0$, this holds if $G(u)\le0$, and, if $G(u)>0$, if and only if $\lambda<\Lambda(u)$. Finally, $\varphi_u(t)\to-\infty$ as $t\to\infty$ when $F(u)>0$.

(iii) Let $\lambda\le\lambda_0^*$ and $F(u)>0$. We claim that $\varphi_u(t)\ge0$ for some $t>0$, with strict inequality if $\lambda<\lambda_0^*$. This is clear if $G(u)\le0$, since then $\varphi_u(t)\ge t^2\big(\frac{1}{2}H(u)-\frac{1}{p}t^{p-2}F(u)\big)>0$ for $t$ small. Let now $u\in V$, so that $\lambda\le\lambda_0^*\le\lambda_0(u)$. We have $\varphi_u(t)=t^q\left(m_u(t)-\frac{\lambda}{q}G(u)\right)$, where $m_u(t)=\frac{1}{2}t^{2-q}H(u)-\frac{1}{p}t^{p-q}F(u)$ is increasing in $(0,t_0(u)]$ and decreasing in $[t_0(u),\infty)$, with
\begin{equation*}
	t_0(u)=\left(\frac{p(2-q)}{2(p-q)}\,\frac{H(u)}{F(u)}\right)^{\frac{1}{p-2}},
\end{equation*}
and a direct computation gives $m_u(t_0(u))=\frac{1}{q}\lambda_0(u)G(u)$. Hence
\begin{equation*}
	\varphi_u(t)\le\frac{1}{q}\left(\lambda_0(u)-\lambda\right)G(u)\,t^q \quad \mbox{for } t>0, \quad \mbox{with equality if, and only if, } t=t_0(u),
\end{equation*}
so the claim holds at $t=t_0(u)$. Note also that if $\varphi_u(t)=\varphi_u'(t)=0$, then $m_u(t)=\frac{\lambda}{q}G(u)$ and $m_u'(t)=0$, so that $t=t_0(u)$ and $\lambda=\lambda_0(u)$. Since $\varphi_u(t)\to0$ as $t\to0$ and $\varphi_u(t)\to-\infty$ as $t\to\infty$, the claim shows that $\varphi_u$ achieves its maximum over $(0,\infty)$ at a local maximizer, which is $t^-(u)$ by (i), and that $\varphi_u(t^-(u))\ge0$, with strict inequality if $\lambda<\lambda_0^*$. Finally, if $\lambda=\lambda_0(u)$, then $\varphi_u\le0=\varphi_u(t_0(u))$, so that $t^-(u)=t_0(u)$ and $\varphi_u(t^-(u))=0$.
\end{proof}

\begin{proof}[Proof of Theorem \ref{thm5}]
By Lemma \ref{lbd0}, $Y=\{F>0\}$, which is open, and (H1''), (H2') and (H3') hold in every path-connected component of $Y$.

(ii) Let $(w_n)\subset\mathcal{N}^-$ be such that $\Phi_\lambda(w_n)\to c:=\inf_{\mathcal{N}^-}\Phi_\lambda$. From $\varphi_{w_n}'(1)=0\ge\varphi_{w_n}''(1)$ we get $H(w_n)=\lambda G(w_n)+F(w_n)$, $(2-q)H(w_n)\le(p-q)F(w_n)$ and
\begin{equation*}
	\Phi_\lambda(w_n)=\left(\frac{1}{2}-\frac{1}{p}\right)H(w_n)-\lambda\left(\frac{1}{q}-\frac{1}{p}\right)G(w_n).
\end{equation*}
Since $|G(w)|\le C\|w\|^q$, $q<2$ and $F(w)\le C\|w\|^p$, $(w_n)$ is bounded and $(F(w_n))$ is bounded away from $0$. If $w_n\rightharpoonup w$, then $F(w)>0$, i.e. $w\in Y$, and since $1=t^-(w_n)$ is a global maximizer of $\varphi_{w_n}$ (Lemma \ref{lbd0}(iii)), $\varphi_w(t)\le\liminf\varphi_{w_n}(t)\le c$ for $t>0$. Thus $w_\lambda:=t^-(w)w\in\mathcal{N}^-$ satisfies $\Phi_\lambda(w_\lambda)=c\le c^-_{W_\lambda}\le\Phi_\lambda(w_\lambda)$, where $W_\lambda$ is the component of $Y$ containing $w$, and $\Phi_\lambda'(w_\lambda)=0$ by Proposition \ref{natural}. If $\lambda<\lambda_0^*$, then $c>0$ by Lemma \ref{lbd0}(iii). Let now $\lambda=\lambda_0^*$, so that $c\ge0$, and let $(u_n)\subset V$ be such that $\|u_n\|=1$, $\lambda_0(u_n)\to\lambda_0^*$ and $u_n\rightharpoonup u$. Since $G$ and $F$ are bounded on the unit sphere, $G(u_n)$ and $F(u_n)$ are bounded away from $0$, so that $u\in V$ and $\lambda_0(u)\le\lim\lambda_0(u_n)=\lambda_0^*$, i.e. $\lambda_0(u)=\lambda$. By Lemma \ref{lbd0}(iii), $\Phi_\lambda(t^-(u)u)=0$, so that $c=0$.

(i) The first assertion follows from Theorem \ref{thmmain}, and, if $\lambda<\lambda_0^*$, the second one from (ii) and Theorem \ref{thmMPG}. {If $\lambda=\lambda_0^*$, the proof of Lemma \ref{lbd0} shows that $t^-(u)>s_u$, so that $\inf_{S\cap Y}t^-\ge\left(\frac{2-q}{(p-q)C}\right)^{\frac1{p-2}}>0$, where $C$ bounds $F$ on $S$. Thus (H5') holds, and Theorem \ref{thmMPG} gives the mountain pass geometry with respect to every component $W$, also when $c_W^-=0$.}
\end{proof}

\subsection{The affine \texorpdfstring{$p$}{p}-Laplacian}\label{subaff}
Since Theorems \ref{thmmain} and \ref{thmMPG} only involve the fibering maps, they also apply to energies that are not coercive. Let $1<p<N$ and let
\begin{equation*}
	\mathcal{E}_p(u)=\gamma_{N,p}\left(\int_{\mathbb{S}^{N-1}}\|\nabla_\xi u\|_p^{-N}\,d\sigma(\xi)\right)^{-\frac{1}{N}}, \quad u\in W^{1,p}_0(\Omega)\setminus\{0\},
\end{equation*}
with $\mathcal{E}_p(0)=0$, be the affine $L^p$ energy of Lutwak, Yang and Zhang \cite{LYZ}. Here $\nabla_\xi u=\nabla u\cdot\xi$ and $\gamma_{N,p}>0$ is a normalizing constant, see \cite{HJM}. The corresponding operator is the affine $p$-Laplacian $\Delta_p^{\mathcal{A}}$ \cite{HJM,LeMo1,LeMo2}. Given $b\in L^\infty(\Omega)$ with $b^+\not\equiv 0$ and $p<r<p^*$, we consider the problem
\begin{equation*}
	\Delta_p^{\mathcal{A}}u=b(x)|u|^{r-2}u \quad \mbox{in } \Omega, \qquad u=0 \quad \mbox{on } \partial\Omega,
\end{equation*}
whose weak solutions are the critical points of
\begin{equation*}
	\Phi_{\mathcal{A}}(u)=\frac{1}{p}\mathcal{E}_p(u)^p-\frac{1}{r}\int_\Omega b|u|^r, \quad u\in W^{1,p}_0(\Omega).
\end{equation*}
The energy $\mathcal{E}_p^p$ is not coercive on $W^{1,p}_0(\Omega)$ \cite{LeMo1}. However, it is $p$-homogeneous, so the fibering maps of $\Phi_{\mathcal{A}}$ have the same structure as for the $p$-Laplacian. Let $Y=\{u\in W^{1,p}_0(\Omega): \int_\Omega b|u|^r>0\}$.

\begin{theorem}\label{thmaff}
Let $1<p<N$, $p<r<p^*$ and $b\in L^\infty(\Omega)$ with $b^+\not\equiv 0$. Then:
\begin{enumerate}[label=(\roman*)]
	\item There exists $c_0>0$ such that, for every path-connected component $W$ of $Y$ and every $e\in W$ with $\Phi_{\mathcal{A}}(e)\le 0$ and $t^-(e)<1$,
	\begin{equation*}
		c_W(e)=c_W^-=\inf_{u\in W}\max_{t>0}\Phi_{\mathcal{A}}(tu)\ge c_0,
	\end{equation*}
	and $\Phi_{\mathcal{A}}$ has the mountain pass geometry with respect to $W$.
	\item There exist a path-connected component $W_0$ of $Y$ and $w\in\mathcal{N}^-\cap W_0$ such that $\Phi_{\mathcal{A}}(w)=c_{W_0}^-=\inf_{\mathcal{N}^-}\Phi_{\mathcal{A}}$ and $\Phi_{\mathcal{A}}'(w)=0$.
\end{enumerate}
\end{theorem}

\begin{remark}
The concave-convex problem $\Delta_p^{\mathcal{A}}u=\lambda a(x)|u|^{q-2}u+b(x)|u|^{r-2}u$, with $1<q<p<r<p^*$ and sign-changing weights, can be treated in the same way. The fibering maps have the same structure as in Theorem \ref{thm5}, with $\frac{1}{2}\int_\Omega|\nabla u|^2$ replaced by $\frac{1}{p}\mathcal{E}_p(u)^p$, and the positivity of $\lambda_0^*$ follows from the affine Sobolev inequality. For positive weights, every fibering map has two critical points, and this problem was studied in \cite{LRS2} by the Nehari manifold method for functionals in $C^1(X\setminus\{0\})$ \cite{LRS1}. Note that the Palais--Smale condition for $\Phi_{\mathcal{A}}$ is not known \cite{LRS2}, which is why we obtain $w$ in Theorem \ref{thmaff}(ii) by minimization rather than through Corollary \ref{existence}.
\end{remark}

Let us prove Theorem \ref{thmaff}. We write $\|u\|=\|\nabla u\|_p$ and $F(u)=\int_\Omega b|u|^r$ for $u\in W^{1,p}_0(\Omega)$, and we use the following facts:
\begin{enumerate}[label=(\alph*)]
	\item $\Phi_{\mathcal{A}}\in C(W^{1,p}_0(\Omega))\cap C^1(W^{1,p}_0(\Omega)\setminus\{0\})$ \cite[Theorem 1]{LeMo2};
	\item the affine Sobolev inequality $\|u\|_{p^*}\le C\mathcal{E}_p(u)$ \cite{LYZ}; in particular $\mathcal{E}_p(u)>0$ for $u\neq 0$;
	\item $\mathcal{E}_p(u)\ge C\|u\|_p^{\frac{N-1}{N}}\|u\|^{\frac{1}{N}}$ for $u\in W^{1,p}_0(\Omega)$ \cite[Theorem 9]{HJM};
	\item sequences that are bounded with respect to $\mathcal{E}_p$ are relatively compact in $L^s(\Omega)$ for $1\le s<p^*$ \cite[Theorem 6.5.3]{Tin};
	\item $\mathcal{E}_p$ is weakly lower semicontinuous on $W^{1,p}_0(\Omega)$ \cite{LeMo1}.
\end{enumerate}

\begin{proof}[Proof of Theorem \ref{thmaff}]
(i) For $u\neq0$ we have $\varphi_u(t)=\frac{t^p}{p}\mathcal{E}_p(u)^p-\frac{t^r}{r}F(u)$, with $\mathcal{E}_p(u)>0$ by (b). Hence $\varphi_u$ has a local maximizer if, and only if, $F(u)>0$. In this case it is the unique critical point of $\varphi_u$, given by $t^-(u)=\left(\mathcal{E}_p(u)^p/F(u)\right)^{\frac{1}{r-p}}$, and
\begin{equation}\label{affenergy}
	\varphi_u(t^-(u))=\left(\frac{1}{p}-\frac{1}{r}\right)\frac{\mathcal{E}_p(u)^{\frac{pr}{r-p}}}{F(u)^{\frac{p}{r-p}}}.
\end{equation}
Thus $Y=\{F>0\}$, which is open, and (H1'') and (H2') hold in every component of $Y$. By (b) and H\"older's inequality, $F(u)\le C\mathcal{E}_p(u)^r$, so \eqref{affenergy} yields $\varphi_u(t^-(u))\ge\left(\frac{1}{p}-\frac{1}{r}\right)C^{-\frac{p}{r-p}}=:c_0>0$ for $u\in Y$, and Theorems \ref{thmmain} and \ref{thmMPG} yield (i).

(ii) Let $(u_n)\subset\mathcal{N}^-$ be such that $\Phi_{\mathcal{A}}(u_n)\to c:=\inf_{\mathcal{N}^-}\Phi_{\mathcal{A}}$. Since $\mathcal{E}_p(u_n)^p=F(u_n)\le C\mathcal{E}_p(u_n)^r$ and $\Phi_{\mathcal{A}}(u_n)=\left(\frac{1}{p}-\frac{1}{r}\right)\mathcal{E}_p(u_n)^p$, $(\mathcal{E}_p(u_n))$ is bounded and bounded away from zero. By (d), up to a subsequence, $u_n\to u$ in $L^s(\Omega)$ for $1\le s<p^*$, so that $F(u)=\lim\mathcal{E}_p(u_n)^p>0$, and since $\|u_n\|_p\to\|u\|_p>0$, (c) shows that $(u_n)$ is bounded in $W^{1,p}_0(\Omega)$, so that $u_n\rightharpoonup u$. With $J(v):=\varphi_v(t^-(v))$, we have $J(u_n)=\Phi_{\mathcal{A}}(u_n)$, since $t^-(u_n)=1$, and (e) and \eqref{affenergy} give $c\le J(u)\le\liminf J(u_n)=c$, so $w:=t^-(u)u$ satisfies $\Phi_{\mathcal{A}}(w)=c\le c_{W_0}^-\le\Phi_{\mathcal{A}}(w)$, where $W_0$ is the component of $Y$ containing $w$. By (a) and Proposition \ref{natural}, $\Phi_{\mathcal{A}}'(w)=0$.
\end{proof}

\subsection{A nonlocal problem}\label{subnl}
We now show that the phenomenon described in Remark \ref{MPGC} also occurs for an elliptic problem, as long as the nonlinearity is allowed to be nonlocal. Let $\lambda_1<\lambda_2$ be the first two eigenvalues of $-\Delta$ in $H^1_0(\Omega)$, and let $\phi_1>0$ be an eigenfunction associated to $\lambda_1$ with $\int_\Omega\phi_1^2=1$. Let $f_1,f_2\in L^2(\Omega)$ be linearly independent, set $\ell_i(u)=\int_\Omega f_iu$ and $a_i=\ell_i(\phi_1)$, $i=1,2$, and assume that
\begin{equation}\label{signs}
	a_1<0<a_2 \quad \mbox{and} \quad a_1+a_2<0.
\end{equation}
Consider
\begin{equation*}
	\Phi_\lambda(u)=\frac{1}{2}\int_\Omega\left(|\nabla u|^2-\lambda u^2\right)-\frac{1}{3}\,\ell_1(u)\ell_2(u)\big(\ell_1(u)+\ell_2(u)\big), \quad u\in H^1_0(\Omega),
\end{equation*}
whose critical points are the weak solutions of the nonlocal problem
\begin{equation*}
	\left\{\begin{array}{ll}
	-\Delta u=\lambda u+\frac{1}{3}\Big[\ell_2(u)\big(2\ell_1(u)+\ell_2(u)\big)f_1+\ell_1(u)\big(\ell_1(u)+2\ell_2(u)\big)f_2\Big] & \mbox{in } \Omega,\\
	u=0 & \mbox{on } \partial\Omega.
	\end{array}\right.
\end{equation*}

\begin{theorem}\label{thm6}
Assume \eqref{signs}. Then there exists $\delta>0$ such that, for every $\lambda\in(\lambda_1,\lambda_1+\delta)$:
\begin{enumerate}[label=(\roman*)]
	\item $\Phi_\lambda(t\phi_1)<0$ for every $t>0$. In particular, there are no $\rho,\alpha>0$ such that $\Phi_\lambda\ge\alpha$ on $\partial B_\rho$, i.e. $\Phi_\lambda$ does not satisfy the classical mountain pass geometry.
	\item The cone $W:=\{u\in H^1_0(\Omega): \ell_1(u)>0,\ \ell_2(u)>0\}$ is a path-connected component of $Y$, and $\Phi_\lambda$ has the mountain pass geometry with respect to $W$. Moreover, $c_W(e)=c_W^-=\inf_{u \in W} \max_{t>0} \Phi_\lambda(tu)>0$ for every $e\in W$ such that $\Phi_\lambda(e)\le 0$ and $t^-(e)<1$, and there exists $w\in\mathcal{N}^-\cap W$ such that $\Phi_\lambda(w)=c_W^-$ and $\Phi_\lambda'(w)=0$.
\end{enumerate}
\end{theorem}

\begin{remark}\label{local}
The nonlinearity in Theorem \ref{thm6} is nonlocal for the following reason. Let $\Phi=\frac{1}{2}H-\frac{1}{p}F$, where $H,F\in C(X)$ are homogeneous of degrees $2$ and $p>2$ respectively, and let $e\in W$ be as in Theorem \ref{thmmain}. Suppose that $\Phi(tv)<0$ for every $t>0$, for some $v\in S$, and that $v$ can be joined to $e/\|e\|$ by a path $\eta$ in $S$ along which $F>0$. Since $F\ge\delta>0$ on $\eta([0,1])$, we have $\Phi(R\eta(s))<0$ for all $s\in[0,1]$ if $R$ is large. Going from $0$ to $Rv$ along the ray of $v$, then from $Rv$ to $Re/\|e\|$ along $R\eta$, and finally back to $e$ along the ray of $e$ (where $\Phi\le 0$ by Lemmas \ref{shape} and \ref{l1}), we obtain a path in $\Gamma$ along which $\Phi\le 0$. Hence $c_W(e)\le 0$, and $\Phi$ cannot have the mountain pass geometry with respect to any cone containing $e$. For local nonlinearities such as $F(u)=\int_\Omega b|u|^p$, such paths are typically available, since mass can be moved between different regions of $\Omega$ by means of functions with disjoint supports. In Theorem \ref{thm6}, as in the planar example of Theorem \ref{thm3.1}, the set of directions where $F>0$ splits into separated pieces: the ray of $\phi_1$ and the cone $W$ lie in different ones, and they are separated by directions where $F\le 0$. Local problems may nevertheless have a disconnected admissible cone when the quadratic part is not even, as for the jumping nonlinearity of Subsection \ref{subfucik}; there the cone is disconnected exactly below the first nontrivial curve of the Fu\v{c}\'{\i}k spectrum (Theorem \ref{thmfucik}).
\end{remark}

Let us prove Theorem \ref{thm6}. We write $\|u\|^2=\int_\Omega|\nabla u|^2$,
\begin{equation*}
	H_\lambda(u)=\int_\Omega\left(|\nabla u|^2-\lambda u^2\right) \quad \mbox{and} \quad F(u)=\ell_1(u)\ell_2(u)\big(\ell_1(u)+\ell_2(u)\big),
\end{equation*}
so that $\Phi_\lambda=\frac{1}{2}H_\lambda-\frac{1}{3}F$, $\varphi_u(t)=\frac{1}{2}H_\lambda(u)t^2-\frac{1}{3}F(u)t^3$ and $\varphi_u'(t)=t\big(H_\lambda(u)-F(u)t\big)$. Thus $\varphi_u$ has a local maximizer if, and only if, $H_\lambda(u)>0$ and $F(u)>0$. In this case it is the unique critical point of $\varphi_u$, given by $t^-(u)=H_\lambda(u)/F(u)$, and
\begin{equation}\label{nlenergy}
	\varphi_u(t^-(u))=\frac{H_\lambda(u)^3}{6F(u)^2}.
\end{equation}
Hence $Y=\{u: H_\lambda(u)>0,\ F(u)>0\}$, and (H1'') and (H2') hold in any path-connected component of $Y$. Note also that $|\ell_i(u)|\le\|f_i\|_2\|u\|_2\le\lambda_1^{-1/2}\|f_i\|_2\|u\|$, so that $|F(u)|\le C\|u\|^3$ for some $C>0$.

\begin{lemma}\label{keyest}
There exist $\delta,c>0$ such that $H_\lambda(u)\ge c\|u\|^2$ for every $u\in\overline{W}$ and $\lambda\in(\lambda_1,\lambda_1+\delta)$.
\end{lemma}

\begin{proof}
Write $u=s\phi_1+w$, with $\int_\Omega w\phi_1=0$, so that $\|u\|^2=\lambda_1s^2+\|w\|^2$, $\int_\Omega w^2\le\lambda_2^{-1}\|w\|^2$ and
\begin{equation*}
	H_\lambda(u)\ge-(\lambda-\lambda_1)s^2+\left(1-\frac{\lambda}{\lambda_2}\right)\|w\|^2.
\end{equation*}
If $u\in\overline{W}$, then $sa_i+\ell_i(w)=\ell_i(u)\ge 0$ for $i=1,2$, and since $|\ell_i(w)|\le\lambda_2^{-1/2}\|f_i\|_2\|w\|$, using $i=1$ if $s\ge0$ and $i=2$ if $s<0$ we get $|s|\le K\|w\|$, where $K=\lambda_2^{-1/2}\max\{\|f_1\|_2/|a_1|,\|f_2\|_2/a_2\}$. Therefore
\begin{equation*}
	H_\lambda(u)\ge\left(1-\frac{\lambda_1}{\lambda_2}-(\lambda-\lambda_1)\Big(\frac{1}{\lambda_2}+K^2\Big)\right)\|w\|^2\ge\frac{1}{2}\left(1-\frac{\lambda_1}{\lambda_2}\right)\|w\|^2
\end{equation*}
whenever $\lambda-\lambda_1<\delta:=\frac{1}{2}\big(1-\frac{\lambda_1}{\lambda_2}\big)\big(\frac{1}{\lambda_2}+K^2\big)^{-1}$. Since $\|u\|^2\le(1+\lambda_1K^2)\|w\|^2$, the conclusion follows.
\end{proof}

\begin{proof}[Proof of Theorem \ref{thm6}]
(i) We have $H_\lambda(\phi_1)=\lambda_1-\lambda<0$ and $F(\phi_1)=a_1a_2(a_1+a_2)>0$ by \eqref{signs}, so that $\Phi_\lambda(t\phi_1)=\frac{1}{2}(\lambda_1-\lambda)t^2-\frac{1}{3}F(\phi_1)t^3<0$ for every $t>0$. In particular, $\Phi_\lambda(\rho\phi_1/\|\phi_1\|)<0$ for every $\rho>0$.

(ii) Let $\delta$ and $c$ be given by Lemma \ref{keyest}. Since $\ell_1,\ell_2$ are linearly independent, $W$ is a nonempty open convex cone, $F>0$ in $W$ and, by Lemma \ref{keyest}, $H_\lambda>0$ in $W$, so that $W\subset Y$. As $\partial W\subset\{F=0\}$ does not meet $Y$, $W$ is a component of $Y$. By \eqref{nlenergy} and Lemma \ref{keyest}, $c_W^-\ge c^3/(6C^2)>0$, i.e. (H3') holds, and (H4') holds: if $S\cap W\ni u_n\to u\in\partial W$, then $t^-(u_n)=H_\lambda(u_n)/F(u_n)\to\infty$, since $H_\lambda(u_n)\ge c$ and $F(u_n)\to0$. Theorems \ref{thmmain} and \ref{thmMPG} yield all the assertions of (ii), except the last one. To prove it, let $(u_n)\subset S\cap W$ be such that $J(u_n)\to c_W^-$, where $J(u):=H_\lambda(u)^3/(6F(u)^2)$, and $u_n\rightharpoonup u$. Since $H_\lambda(u_n)\ge c$ and $(J(u_n))$ is bounded, $F(u)=\lim F(u_n)>0$, which forces $\ell_i(u)>0$, since $\ell_i(u)\ge0$; thus $u\in W$ and $c_W^-\le J(u)\le\liminf J(u_n)=c_W^-$, since $J$ is homogeneous of degree zero. Thus $w:=t^-(u)u\in\mathcal{N}^-\cap W$ satisfies $\Phi_\lambda(w)=c_W^-$, and $\Phi_\lambda'(w)=0$ by Proposition \ref{natural}.
\end{proof}

\subsection{Finite-rank problems and multiplicity}\label{subfr}
The finite-rank structure used in Subsection \ref{subnl} also produces problems with many components, and hence, by Corollary \ref{multiplicity}, with many critical points. Let $m\ge2$, $\lambda<\lambda_1$, and let $\phi_1,\dots,\phi_m$ be $L^2(\Omega)$-orthonormal Dirichlet eigenfunctions of $-\Delta$, associated with the eigenvalues $\lambda_1\le\dots\le\lambda_m$. Set
\begin{equation*}
	Q(u)=\int_\Omega\left(|\nabla u|^2-\lambda u^2\right), \quad c_i(u)=\int_\Omega u\phi_i, \quad q_i=\lambda_i-\lambda>0,
\end{equation*}
and consider
\begin{equation*}
	\Phi_m(u)=\frac{1}{2}Q(u)-\frac{1}{2m}\prod_{i=1}^mc_i(u)^2, \quad u\in H^1_0(\Omega),
\end{equation*}
whose critical points are the weak solutions of the nonlocal problem
\begin{equation*}
	-\Delta u-\lambda u=\frac{1}{m}\sum_{i=1}^m c_i(u)\Big(\prod_{j\neq i}c_j(u)^2\Big)\phi_i \quad \mbox{in } \Omega, \qquad u=0 \quad \mbox{on } \partial\Omega.
\end{equation*}
For $\sigma=(\sigma_1,\dots,\sigma_m)\in\{-1,1\}^m$, let $C_\sigma=\{u\in H^1_0(\Omega): \sigma_ic_i(u)>0,\ i=1,\dots,m\}$.

\begin{theorem}\label{thmfr}
The sets $C_\sigma$, $\sigma\in\{-1,1\}^m$, are the $2^m$ path-connected components of $Y$, and the assumptions of Corollary \ref{multiplicity} hold in each of them. The nonzero critical points of $\Phi_m$ are exactly
\begin{equation*}
	u_\sigma=\sum_{i=1}^m\sigma_i\left(\frac{P}{mq_i}\right)^{\frac{1}{2}}\phi_i, \quad \sigma\in\{-1,1\}^m, \quad \mbox{where} \quad P=\Big(m^m\prod_{i=1}^mq_i\Big)^{\frac{1}{m-1}}.
\end{equation*}
Moreover, $u_\sigma\in C_\sigma$ and
\begin{equation*}
	\Phi_m(u_\sigma)=c^-_{C_\sigma}=\inf_{\mathcal{N}}\Phi_m=\frac{m-1}{2m}P \quad \mbox{for every } \sigma.
\end{equation*}
\end{theorem}

The levels coincide because $\Phi_m$ is invariant under independent changes of sign of the coefficients $c_i(u)$. Breaking this invariance yields distinct levels, while Corollary \ref{multiplicity} still provides one critical point in each component.

\begin{proposition}\label{propfr}
Let $m=2$, $0<\delta<\varepsilon$ with $\varepsilon+\delta<1$, and
\begin{equation*}
	\tilde\Phi(u)=\frac{1}{2}Q(u)-\frac{1}{4}F_{\varepsilon,\delta}(c_1(u),c_2(u)), \quad F_{\varepsilon,\delta}(a,b)=a^2b^2+\varepsilon a|a|b^2+\delta a^2b|b|.
\end{equation*}
Then $\tilde\Phi\in C^1(H^1_0(\Omega))$, the components of $Y$ are again the four sets $C_\sigma$, and the assumptions of Corollary \ref{multiplicity} hold in each of them. The only critical point of $\tilde\Phi$ in $C_\sigma$ is
\begin{equation*}
	\tilde u_\sigma=\sigma_1\sqrt{\frac{2q_2}{\kappa_\sigma}}\,\phi_1+\sigma_2\sqrt{\frac{2q_1}{\kappa_\sigma}}\,\phi_2, \quad \mbox{where} \quad \kappa_\sigma=1+\varepsilon\sigma_1+\delta\sigma_2,
\end{equation*}
and $\tilde\Phi(\tilde u_\sigma)=c^-_{C_\sigma}=q_1q_2/\kappa_\sigma$. In particular, $c^-_{C_{++}}<c^-_{C_{+-}}<c^-_{C_{-+}}<c^-_{C_{--}}$.
\end{proposition}

\begin{remark}
What matters in Corollary \ref{multiplicity} is the number of components of $Y$, not the number of modes in the equation: $m$ modes produce $2^m$ components, indexed by the sign patterns of $(c_1(u),\dots,c_m(u))$. Proposition \ref{propfr} shows that the multiplicity does not depend on the levels being equal. Since the critical points of $\Phi_m$ and $\tilde\Phi$ are explicit, the interest of these examples lies in the minimax characterization of each of them, rather than in their existence.
\end{remark}

\begin{proof}[Proof of Theorem \ref{thmfr}]
We use the norm $\|u\|=Q(u)^{1/2}$, which is equivalent to the usual one since $\lambda<\lambda_1$, and write $P(u)=\prod_{i=1}^mc_i(u)^2$. Since $\varphi_u(t)=\frac{t^2}{2}Q(u)-\frac{t^{2m}}{2m}P(u)$, $\varphi_u$ has a local maximizer if, and only if, $P(u)>0$, and it is then its unique critical point, with $\varphi_u(t^-(u))=\frac{m-1}{2m}Q(u)^{\frac{m}{m-1}}P(u)^{-\frac{1}{m-1}}$. Thus $Y=\{P>0\}$ is the union of the open convex sets $C_\sigma$, which are its components since $\partial C_\sigma\subset\{P=0\}$, (H1'') and (H2') hold in each $C_\sigma$, and $c^-_{C_\sigma}>0$, since $P$ is bounded on the unit sphere. Palais--Smale sequences are bounded, since $\Phi_m(u)-\frac{1}{2m}\Phi_m'(u)u=\frac{m-1}{2m}Q(u)$, and if $u_n\rightharpoonup u$, then $c_i(u_n)\to c_i(u)$, so the nonlinear part of $\Phi_m'(u_n)$, which lies in $\mathrm{span}\{\phi_1,\dots,\phi_m\}$, converges, and $\Phi_m'(u_n)\to0$ yields $u_n\to u$. Hence Corollary \ref{multiplicity} applies in each $C_\sigma$.

Let $u\neq0$ be a critical point of $\Phi_m$. Testing $\Phi_m'(u)=0$ with the component $z$ of $u$ that is $L^2$-orthogonal to $\mathrm{span}\{\phi_1,\dots,\phi_m\}$ gives $Q(z)=0$, so that $u=\sum_{i=1}^ma_i\phi_i$, and testing with $\phi_i$ gives $q_ia_i=\frac{1}{m}a_i\prod_{j\neq i}a_j^2$. If some $a_i$ vanished, all the right-hand sides would vanish and $u=0$. Hence all the $a_i$ are nonzero, so $a_i^2=\frac{P}{mq_i}$ with $P=\prod_ja_j^2$, and taking the product over $i$ gives the value of $P$ in the statement. Thus the nonzero critical points are exactly the points $u_\sigma$, and $u_\sigma$ is the only one in $C_\sigma$, so that Corollary \ref{multiplicity} gives $c^-_{C_\sigma}=\Phi_m(u_\sigma)=\frac{m-1}{2m}P$, since $Q(u_\sigma)=\sum_iq_ia_i^2=P$. Finally, $\Phi_m'(u)u=Q(u)-P(u)$, so $\mathcal{N}\subset Y$, and $\inf_{\mathcal{N}}\Phi_m=\min_\sigma c^-_{C_\sigma}=\frac{m-1}{2m}P$.
\end{proof}

\begin{proof}[Proof of Proposition \ref{propfr}]
Since $t\mapsto t|t|$ is of class $C^1$, so is $\tilde\Phi$. In $C_\sigma$ we have $F_{\varepsilon,\delta}(c_1(u),c_2(u))=\kappa_\sigma c_1(u)^2c_2(u)^2$, with $\kappa_\sigma>0$ since $\varepsilon+\delta<1$, while $F_{\varepsilon,\delta}(c_1(u),c_2(u))=0$ if $c_1(u)c_2(u)=0$. Hence the proof of Theorem \ref{thmfr} applies, with the nonlinear term multiplied by $\kappa_\sigma$ in $C_\sigma$ and $\tilde\Phi(u)-\frac{1}{4}\tilde\Phi'(u)u=\frac{1}{4}Q(u)$. A critical point $u\in C_\sigma$ is of the form $u=a\phi_1+b\phi_2$, with $q_1a=\frac{\kappa_\sigma}{2}ab^2$ and $q_2b=\frac{\kappa_\sigma}{2}a^2b$, so it is $\tilde u_\sigma$, and $c^-_{C_\sigma}=\tilde\Phi(\tilde u_\sigma)=q_1q_2/\kappa_\sigma$ by Corollary \ref{multiplicity}. Finally, $0<\delta<\varepsilon$ gives $\kappa_{++}>\kappa_{+-}>\kappa_{-+}>\kappa_{--}$.
\end{proof}

\subsection{A local problem related to the Fu\v{c}\'{\i}k spectrum}\label{subfucik}
We now give a local problem whose admissible cone is disconnected exactly below the first nontrivial curve of the Fu\v{c}\'{\i}k spectrum. Let $\Omega\subset\mathbb{R}^N$ be a bounded smooth domain, let $\lambda_1$ be the first Dirichlet eigenvalue of $-\Delta$, and let $\phi_1>0$ be an associated eigenfunction with $\|\phi_1\|_2=1$. Recall that the Fu\v{c}\'{\i}k spectrum $\Sigma$ is the set of $(\alpha,\beta)\in\mathbb{R}^2$ for which $-\Delta u=\alpha u^+-\beta u^-$ has a nontrivial solution in $H^1_0(\Omega)$. For $s\ge 0$, let
\begin{equation*}
	J_s(u)=\int_\Omega|\nabla u|^2-s\int_\Omega(u^+)^2, \qquad S_2=\{u\in H^1_0(\Omega): \|u\|_2=1\},
\end{equation*}
and
\begin{equation}\label{cs}
	c(s)=\inf_{\gamma\in\Gamma_s}\max_{t\in[-1,1]}J_s(\gamma(t)), \qquad \Gamma_s=\{\gamma\in C([-1,1],S_2): \gamma(-1)=-\phi_1,\ \gamma(1)=\phi_1\}.
\end{equation}
Then $(s+c(s),c(s))$ belongs to the first nontrivial curve of $\Sigma$, and $c(s)>\lambda_1$ \cite{dFG,DP}. Fix $s>0$, $\beta>\lambda_1$, $\alpha=\beta+s$, $1<q<2$ and $\mu_\pm>0$, and consider
\begin{equation}\label{fucikpb}
	-\Delta u=\alpha u^+-\beta u^--\mu_+(u^+)^{q-1}+\mu_-(u^-)^{q-1} \quad \mbox{in } \Omega, \qquad u=0 \quad \mbox{on } \partial\Omega,
\end{equation}
whose weak solutions are the critical points of $\Phi=\frac{1}{2}Q_{\alpha,\beta}+\frac{1}{q}B\in C^1(H^1_0(\Omega))$, where
\begin{equation*}
	Q_{\alpha,\beta}(u)=\int_\Omega|\nabla u|^2-\alpha\int_\Omega(u^+)^2-\beta\int_\Omega(u^-)^2, \qquad B(u)=\mu_+\int_\Omega(u^+)^q+\mu_-\int_\Omega(u^-)^q.
\end{equation*}
Since $\alpha\neq\beta$, the functional $\Phi$ is not even. Let
\begin{equation*}
	d_+=\inf\{\Phi(u): u\in\mathcal{N}^-,\ u\ge0\} \quad \mbox{and} \quad d_-=\inf\{\Phi(u): u\in\mathcal{N}^-,\ u\le0\}.
\end{equation*}
For $u\ge0$, \eqref{fucikpb} reduces to the problem with strong absorption $-\Delta u+\mu_+u^{q-1}=\alpha u$ studied in \cite{DHI}, and $d_+$ is the energy of its ground states. Similarly, $d_-$ is the energy of the ground states of $-\Delta w+\mu_-w^{q-1}=\beta w$, where $w=-u\ge0$. Both levels are positive and achieved for every $\alpha,\beta>\lambda_1$ \cite[Theorem 2.1]{DHI}.

\begin{theorem}\label{thmfucik}
Under the above assumptions, $Y=\{u\in H^1_0(\Omega): Q_{\alpha,\beta}(u)<0\}$, and (H1''), (H2') and (H3') hold in every path-connected component of $Y$. Moreover:
\begin{enumerate}[label=(\roman*)]
	\item If $\beta\le c(s)$, then $Y$ has exactly two path-connected components $C_+$ and $C_-$, which contain, respectively, the nonnegative and the nonpositive functions of $Y$, and
	\begin{equation*}
		c(e)=c^-_{C_\pm}=d_\pm \quad \mbox{for every } e\in C_\pm \mbox{ such that } \Phi(e)\le0 \mbox{ and } t^-(e)<1.
	\end{equation*}
	\item If $\beta>c(s)$, then $Y$ is path-connected, and
	\begin{equation*}
		c(e)=c^-_Y=\min\{d_+,d_-\} \quad \mbox{for every } e\in Y \mbox{ such that } \Phi(e)\le0 \mbox{ and } t^-(e)<1.
	\end{equation*}
	\item If $\beta<c(s)$, then the assumptions of Corollary \ref{multiplicity} hold in $C_\pm$, and \eqref{fucikpb} has two nontrivial solutions $u_\pm\in\mathcal{N}^-\cap C_\pm$ such that $u_+\ge0$, $u_-\le0$ and $\Phi(u_\pm)=c^-_{C_\pm}=d_\pm$.
\end{enumerate}
\end{theorem}

\begin{remark}\label{remfucik}
Theorem \ref{thmfucik} shows that the first nontrivial curve of $\Sigma$ is a threshold for the minimax structure of $\Phi$. If $\beta\le c(s)$, the mountain pass level between $0$ and $e$ is $d_+$ or $d_-$, according to the component of $Y$ that contains $e$, while if $\beta>c(s)$ it equals $\min\{d_+,d_-\}$ for every $e$. Hence, whenever $d_+\neq d_-$, which is the case, for instance, if $\mu_+\le\mu_-$ (Proposition \ref{propfucik}), the mountain pass level depends on the endpoint $e$ if, and only if, $\beta\le c(s)$. The disconnection of $Y$ comes from the minimax characterization \eqref{cs} of $c(s)$ and the identity $Q_{\alpha,\beta}(u)=J_s(u)-\beta\|u\|_2^2$. Parts (i) and (ii) also hold at $\beta=c(s)$, where $(\alpha,\beta)\in\Sigma$, since Theorem \ref{thmmain} requires neither the Palais--Smale condition nor $(\alpha,\beta)\notin\Sigma$.

The existence of the solutions $u_\pm$ is not new: they are ground states of the two problems with strong absorption mentioned above, which have nonnegative ground states for all $\alpha,\beta>\lambda_1$ \cite[Theorem 2.1]{DHI}; see also \cite{Por}, where existence is obtained by a variant of the mountain pass theorem. By \cite[Proposition 5.1 and Theorem 5.1]{DHI}, such ground states may have compact support, so that $u_+>0$ cannot be expected in general. The contribution of Theorem \ref{thmfucik} is the structure of $Y$ and the resulting minimax characterization of $d_\pm$, by paths that are not required to stay in $C_\pm$, nor to have a sign.
\end{remark}

\begin{proof}[Proof of Theorem \ref{thmfucik}]
\emph{Step 1: the fiber maps.} For $u\neq0$ we have $\varphi_u(t)=\frac{t^2}{2}Q_{\alpha,\beta}(u)+\frac{t^q}{q}B(u)$ and $B(u)>0$. Since $q<2$, $\varphi_u$ has a local maximizer if, and only if, $Q_{\alpha,\beta}(u)<0$. In this case it is the unique critical point of $\varphi_u$, it is a global maximizer, and
\begin{equation*}
	t^-(u)=\left(\frac{B(u)}{-Q_{\alpha,\beta}(u)}\right)^{\frac{1}{2-q}}, \qquad \varphi_u(t^-(u))=\Big(\frac{1}{q}-\frac{1}{2}\Big)B(u)^{\frac{2}{2-q}}\big(-Q_{\alpha,\beta}(u)\big)^{-\frac{q}{2-q}}>0.
\end{equation*}
Thus $Y=\{Q_{\alpha,\beta}<0\}$, and (H1'') and (H2') hold in every component of $Y$. The function $E(u):=\varphi_u(t^-(u))$, $u\in Y$, is homogeneous of degree zero and increasing in $B(u)$ and in $Q_{\alpha,\beta}(u)$, $c_W^-=\inf_WE$ for every component $W$ of $Y$, $d_+=\inf\{E(u): u\in Y,\ u\ge0\}$ and $d_-=\inf\{E(u): u\in Y,\ u\le0\}$.

\emph{Step 2: the components of $Y$.} Since $\alpha=\beta+s$, we have $Q_{\alpha,\beta}(u)=J_s(u)-\beta\|u\|_2^2$, so that $Y\cap S_2=\{u\in S_2: J_s(u)<\beta\}$ contains $\pm\phi_1$. Let $C_\pm$ be the component of $Y$ containing $\pm\phi_1$. Recall that $Q_{\alpha,\beta}(u)=Q_{\alpha,\beta}(u^+)+Q_{\alpha,\beta}(-u^-)$, $B(u)=B(u^+)+B(-u^-)$, and that $Q_{\alpha,\beta}$ coincides with $q_\alpha(v)=\|\nabla v\|_2^2-\alpha\|v\|_2^2$ on nonnegative functions and with $q_\beta$ on nonpositive ones. If $v\in Y$ is nonnegative, then, since $\int_\Omega\nabla v\cdot\nabla\phi_1=\lambda_1\int_\Omega v\phi_1>0$, the functions $z_t=(1-t)v+t\phi_1\ge0$ satisfy
\begin{equation*}
	Q_{\alpha,\beta}(z_t)=(1-t)^2q_\alpha(v)+t^2(\lambda_1-\alpha)+2t(1-t)(\lambda_1-\alpha)\int_\Omega v\phi_1<0, \qquad t\in[0,1],
\end{equation*}
so that $v\in C_+$; similarly, the nonpositive functions of $Y$ lie in $C_-$. If $u\in Y$ changes sign and $Q_{\alpha,\beta}(-u^-)<0$, then $Q_{\alpha,\beta}(su^+-u^-)=s^2Q_{\alpha,\beta}(u^+)+Q_{\alpha,\beta}(-u^-)\le\max\{Q_{\alpha,\beta}(-u^-),Q_{\alpha,\beta}(u)\}<0$ for $s\in[0,1]$, so that $u\in C_-$. Similarly, $u\in C_+$ if $Q_{\alpha,\beta}(u^+)<0$, and one of these cases occurs, since $Q_{\alpha,\beta}(u^+)+Q_{\alpha,\beta}(-u^-)=Q_{\alpha,\beta}(u)<0$. Hence $Y=C_+\cup C_-$. If $\beta\le c(s)$, then $C_+\neq C_-$, since normalizing in $L^2$ a path in $Y$ from $-\phi_1$ to $\phi_1$ would give $\gamma\in\Gamma_s$ with $J_s<\beta$ on the compact set $\gamma([-1,1])$, so that $c(s)<\beta$; consequently, $Q_{\alpha,\beta}(-u^-)\ge0$ for $u\in C_+$ and $Q_{\alpha,\beta}(u^+)\ge0$ for $u\in C_-$. If $\beta>c(s)$, then $J_s<\beta$ on $\gamma([-1,1])$ for some $\gamma\in\Gamma_s$, i.e. $\gamma([-1,1])\subset Y$, so that $C_+=C_-=Y$.

\emph{Step 3: (H3').} If $u\in Y$ and $\|\nabla u\|_2=1$, then $1<\alpha\|u^+\|_2^2+\beta\|u^-\|_2^2\le\alpha\|u\|_2^2$, so that, by interpolation between $L^q$ and $L^r$, for some $r\in(2,2^*)$, and the Sobolev embedding, $\|u\|_q\ge\eta>0$ and $B(u)\ge\min\{\mu_+,\mu_-\}\eta^q=:b_0$. Since also $-Q_{\alpha,\beta}(u)\le\alpha/\lambda_1$, Step 1 gives $E(u)\ge\big(\frac{1}{q}-\frac{1}{2}\big)b_0^{\frac{2}{2-q}}(\alpha/\lambda_1)^{-\frac{q}{2-q}}$, so that, by homogeneity, $c_W^->0$ for every component $W$ of $Y$.

\emph{Step 4: the levels.} Let $u\in Y$ be such that $u^+\neq0\neq u^-$. If $Q_{\alpha,\beta}(-u^-)\ge0$, then $Q_{\alpha,\beta}(u^+)\le Q_{\alpha,\beta}(u)<0$ and $B(u^+)<B(u)$, so that $u^+\in Y$ and $E(u)>E(u^+)\ge d_+$. Similarly, if $Q_{\alpha,\beta}(u^+)\ge0$, then $E(u)>E(-u^-)\ge d_-$. Assume now that $Q_{\alpha,\beta}(u^+)<0$ and $Q_{\alpha,\beta}(-u^-)<0$, and let $a=\frac{2}{2-q}$, $\theta=\frac{a-1}{a}\in(0,1)$, $B_1=B(u^+)$, $B_2=B(-u^-)$, $D_1=-Q_{\alpha,\beta}(u^+)$ and $D_2=-Q_{\alpha,\beta}(-u^-)$, so that $E(u)=\kappa(B_1+B_2)^a(D_1+D_2)^{1-a}$, where $\kappa=\frac{1}{q}-\frac{1}{2}$. If $\kappa\ell=\min\{E(u^+),E(-u^-)\}$, then $B_i\ge\ell^{1/a}D_i^\theta$ for $i=1,2$, so that
\begin{equation*}
	B_1+B_2\ge\ell^{1/a}\left(D_1^\theta+D_2^\theta\right)\ge\ell^{1/a}(D_1+D_2)^\theta,
\end{equation*}
i.e. $E(u)\ge\min\{E(u^+),E(-u^-)\}\ge\min\{d_+,d_-\}$. Since the nonnegative and the nonpositive functions of $Y$ belong to $Y$, we conclude that $\inf_YE=\min\{d_+,d_-\}$. If $\beta\le c(s)$, then by Step 2 every $u\in C_+$ with $u^-\neq0$ satisfies $u^+\neq0$ and $Q_{\alpha,\beta}(-u^-)\ge0$, so that $E(u)>d_+$; since the nonnegative functions of $Y$ lie in $C_+$, we get $c^-_{C_+}=\inf_{C_+}E=d_+$, and similarly $c^-_{C_-}=d_-$. Assertions (i) and (ii) now follow from Steps 1--4 and Theorem \ref{thmmain}.

\emph{Step 5: if $\beta<c(s)$, then $(\alpha,\beta)\notin\Sigma$.} This holds because $(s+c(s),c(s))$ is the first nontrivial point of $\Sigma$ on the line $\{\alpha-\beta=s\}$ \cite{dFG}.

\emph{Step 6: if $\beta<c(s)$, then $\Phi$ satisfies the Palais--Smale condition.} If $(u_n)$ is a Palais--Smale sequence with $\rho_n=\|u_n\|\to\infty$, then, up to a subsequence, $v_n=u_n/\rho_n\rightharpoonup v$ in $H^1_0(\Omega)$ and $v_n\to v$ in $L^2(\Omega)$, and since $\rho_n^{q-2}\to0$, $v$ solves $-\Delta v=\alpha v^+-\beta v^-$. Testing $\Phi'(u_n)$ with $u_n/\rho_n^2$ gives $1\le\alpha\|v^+\|_2^2+\beta\|v^-\|_2^2$, so $v\neq0$, which contradicts Step 5. Hence Palais--Smale sequences are bounded, and they have convergent subsequences by the compact embeddings into $L^2(\Omega)$ and $L^q(\Omega)$.

\emph{Step 7: proof of (iii).} Let $\beta<c(s)$. By Steps 1, 3 and 6, the assumptions of Corollary \ref{multiplicity} hold in $C_+$ and $C_-$, which yields $u_\pm\in\mathcal{N}^-\cap C_\pm$ such that $\Phi(u_\pm)=c^-_{C_\pm}=d_\pm$. If $u_+^-\neq0$, then Step 4 gives $\Phi(u_+)=E(u_+)>d_+$, a contradiction. Hence $u_+\ge0$, and similarly $u_-\le0$.
\end{proof}

We now compare $d_+$ and $d_-$. Let $\kappa=\frac{1}{q}-\frac{1}{2}$ and $b=\frac{q}{2-q}$, and for $\gamma>\lambda_1$ let
\begin{equation*}
	m(\gamma)=\kappa\inf\left\{\Big(\int_\Omega|v|^q\Big)^{\frac{2}{2-q}}\big(\gamma-\|\nabla v\|_2^2\big)^{-b}: v\in S_2,\ \|\nabla v\|_2^2<\gamma\right\}.
\end{equation*}
Then $\mu^{\frac{2}{2-q}}m(\gamma)$ is the ground state level of $-\Delta u+\mu|u|^{q-2}u=\gamma u$ considered in \cite{DHI}. In particular, $m(\gamma)$ is achieved \cite[Lemma 2.2]{DHI}, $m$ is strictly decreasing and differentiable in $(\lambda_1,\infty)$, with $m'(\gamma)=-\frac{1}{2}\int_\Omega u_\gamma^2$, where $u_\gamma$ is a ground state of $-\Delta u+|u|^{q-2}u=\gamma u$ \cite[Lemma 4.1 and Corollary 4.3]{DHI}, and $m(\gamma)\to\infty$ as $\gamma\to\lambda_1^+$ \cite[Proposition 2.1]{DHI}. The next result makes the monotonicity of $m$ quantitative and compares the two levels.

\begin{proposition}\label{propfucik}
Under the above assumptions, for every $\beta>\lambda_1$:
\begin{enumerate}[label=(\roman*)]
	\item $d_+=\mu_+^{\frac{2}{2-q}}m(\alpha)$ and $d_-=\mu_-^{\frac{2}{2-q}}m(\beta)$.
	\item The map $\gamma\mapsto(\gamma-\lambda_1)^b\,m(\gamma)$ is nonincreasing in $(\lambda_1,\infty)$, and
	\begin{equation*}
		m(\alpha)<\Big(\frac{\beta-\lambda_1}{\alpha-\lambda_1}\Big)^b m(\beta).
	\end{equation*}
	\item Let $\tau=\big(m(\beta)/m(\alpha)\big)^{\frac{2-q}{2}}$. Then $d_+<d_-$, $d_+=d_-$ or $d_+>d_-$, according as $\mu_+/\mu_-<\tau$, $\mu_+/\mu_-=\tau$ or $\mu_+/\mu_->\tau$. Moreover,
	\begin{equation*}
		\Big(\frac{\alpha-\lambda_1}{\beta-\lambda_1}\Big)^{\frac{q}{2}}<\tau\le\frac{\int_\Omega\phi_1^q}{\eta_\alpha}\Big(\frac{\alpha-\lambda_1}{\beta-\lambda_1}\Big)^{\frac{q}{2}}, \quad \mbox{where} \quad \eta_\alpha=\inf\Big\{\int_\Omega|v|^q: v\in S_2,\ \|\nabla v\|_2^2\le\alpha\Big\}>0.
	\end{equation*}
	In particular, $d_+<d_-$ whenever $\mu_+\le\mu_-$.
\end{enumerate}
\end{proposition}

\begin{proof}
For $\gamma>\lambda_1$ and $v\in S_2$ with $\|\nabla v\|_2^2<\gamma$, let $f_\gamma(v)=\kappa\big(\int_\Omega|v|^q\big)^{\frac{2}{2-q}}\big(\gamma-\|\nabla v\|_2^2\big)^{-b}$, so that $m(\gamma)=\inf f_\gamma$, and the infimum may be taken over nonnegative functions.

(i) By Step 1 of the proof of Theorem \ref{thmfucik}, $d_+=\inf\{E(u): u\in Y,\ u\ge0\}$. For $u\ge0$, we have $u\in Y$ if, and only if, $\|\nabla u\|_2^2<\alpha\|u\|_2^2$, and $E(u)=\mu_+^{\frac{2}{2-q}}f_\alpha(u/\|u\|_2)$, since $E$ is homogeneous of degree zero. Thus $d_+=\mu_+^{\frac{2}{2-q}}m(\alpha)$, and similarly $d_-=\mu_-^{\frac{2}{2-q}}m(\beta)$.

(ii) Let $\lambda_1<\gamma_1<\gamma_2$ and $v\in S_2$ with $a:=\|\nabla v\|_2^2<\gamma_1$. Since $a\ge\lambda_1$ and $a\mapsto(\gamma_1-a)/(\gamma_2-a)$ is decreasing,
\begin{equation*}
	f_{\gamma_2}(v)=\Big(\frac{\gamma_1-a}{\gamma_2-a}\Big)^bf_{\gamma_1}(v)\le\Big(\frac{\gamma_1-\lambda_1}{\gamma_2-\lambda_1}\Big)^bf_{\gamma_1}(v),
\end{equation*}
with strict inequality if $a>\lambda_1$, and taking the infimum over $v$ gives the monotonicity. Let $v_\beta\ge0$ be such that $f_\beta(v_\beta)=m(\beta)$ \cite[Lemma 2.2]{DHI}. A positive multiple of $v_\beta$ solves $-\Delta w+w^{q-1}=\beta w$ \cite[Lemma 2.1]{DHI}, so $v_\beta$ is not a multiple of $\phi_1$, since $\phi_1^{q-2}$ is not constant. Hence $\|\nabla v_\beta\|_2^2>\lambda_1$, and $m(\alpha)\le f_\alpha(v_\beta)<\big(\frac{\beta-\lambda_1}{\alpha-\lambda_1}\big)^bm(\beta)$.

(iii) By (i), $d_+/d_-=\big(\mu_+/(\tau\mu_-)\big)^{\frac{2}{2-q}}$, which gives the first assertion, and the lower bound for $\tau$ follows from (ii). For the upper bound, $m(\beta)\le f_\beta(\phi_1)=\kappa\big(\int_\Omega\phi_1^q\big)^{\frac{2}{2-q}}(\beta-\lambda_1)^{-b}$, while $m(\alpha)\ge\kappa\,\eta_\alpha^{\frac{2}{2-q}}(\alpha-\lambda_1)^{-b}$, since $\int_\Omega|v|^q\ge\eta_\alpha$ and $\alpha-\|\nabla v\|_2^2\le\alpha-\lambda_1$; and $\eta_\alpha>0$ by the interpolation argument of Step 3 of the proof of Theorem \ref{thmfucik}.
\end{proof}

\begin{remark}
By Step 4 of the proof of Theorem \ref{thmfucik}, the ground state level $\inf_{\mathcal{N}^-}\Phi=\inf_YE$ equals $\min\{d_+,d_-\}$ for every $\beta>\lambda_1$. Thus, if $\beta>c(s)$, the mountain pass level coincides with the ground state level for every endpoint $e$, whereas if $\beta\le c(s)$ this happens only for the endpoints in the component of $Y$ with the smaller level. For $\beta<c(s)$ and $\mu_+/\mu_-=\tau$, the solutions $u_+\ge0$ and $u_-\le0$ are ground states at the same level, and they are distinguished only by the component of $Y$ in which they lie (see Remark \ref{multrem}). They are not related by the map $u\mapsto-u$, which is not a symmetry of $\Phi$: in fact, $u_-\neq-u_+$ for all values of the parameters. Indeed, if $w=u_+=-u_-$, then $w\ge0$, $w\not\equiv0$, and subtracting the equations $-\Delta w+\mu_+w^{q-1}=\alpha w$ and $-\Delta w+\mu_-w^{q-1}=\beta w$ gives $sw=(\mu_+-\mu_-)w^{q-1}$, which is impossible on $\{w>0\}$ if $\mu_+\le\mu_-$, and otherwise forces $w$ to take only the values $0$ and $\big((\mu_+-\mu_-)/s\big)^{\frac{1}{2-q}}$, so that $\nabla w=0$ a.e. and $w\equiv0$, a contradiction. The same argument shows that the two equations have no common nontrivial solution. For $s=0$ and $\mu_+=\mu_-$, on the other hand, $\Phi$ is even and one may take $u_-=-u_+$. Finally, the estimates in the proof of Proposition \ref{propfucik} give
\begin{equation*}
	\kappa\,\eta_\gamma^{\frac{2}{2-q}}(\gamma-\lambda_1)^{-b}\le m(\gamma)\le\kappa\Big(\int_\Omega\phi_1^q\Big)^{\frac{2}{2-q}}(\gamma-\lambda_1)^{-b}, \qquad \gamma>\lambda_1,
\end{equation*}
where $\eta_\gamma$ is defined as $\eta_\alpha$, with $\alpha$ replaced by $\gamma$, and is nonincreasing in $\gamma$. By the asymptotic profile of the ground states obtained in \cite[Proposition 4.1]{DHI}, the upper bound is sharp as $\gamma\to\lambda_1^+$. In particular, for fixed $s$ and $\mu_\pm$, we have $d_-\to\infty$ as $\beta\to\lambda_1^+$, while $d_+\le\kappa\big(\mu_+\int_\Omega\phi_1^q\big)^{\frac{2}{2-q}}s^{-b}$ remains bounded.
\end{remark}

\subsection{Prescribed scalar curvature}\label{subcurv}
Let $(M,g)$ be a closed Riemannian manifold of dimension $n\ge 3$ with constant scalar curvature $R_g<0$, and let $K\in C(M)$ with $K^+\not\equiv 0$. If $u>0$ solves
\begin{equation}\label{yamabe}
	-a_n\Delta_g u+R_gu=Ku^{\frac{n+2}{n-2}} \quad \mbox{in } M, \qquad a_n=\frac{4(n-1)}{n-2},
\end{equation}
then the conformal metric $u^{\frac{4}{n-2}}g$ has scalar curvature $K$ \cite{KW}. The associated functional is
\begin{equation*}
	\Phi(u)=\frac{1}{2}Q(u)-\frac{1}{2^*}F(u), \quad Q(u)=\int_M\left(a_n|\nabla u|^2+R_gu^2\right), \quad F(u)=\int_MK|u|^{2^*},
\end{equation*}
for $u\in H^1(M)$, where $2^*=\frac{2n}{n-2}$. The quadratic form $Q$ is indefinite, since $Q(1)=R_g|M|<0$, and $F$ involves the critical Sobolev exponent, so $\Phi$ does not satisfy the Palais--Smale condition in general. Since Theorem \ref{thmmain} requires no compactness, it applies to $\Phi$ as follows.

\begin{proposition}\label{curv}
{
Let $D=\{u\in H^1(M): Q(u)>0,\ F(u)>0\}$. Then $D$ is a nonempty open cone and, for every path-connected component $W$ of $D$ and every $e\in W$ with $\Phi(e)\le 0$ and $t^-(e)<1$,

\begin{equation*}
	c_W(e)=c_W^-=\inf_{u\in W}\max_{t>0}\Phi(tu).
\end{equation*}
Moreover, if $\int_MK\ge 0$, then $\Phi(t)<0$ for every constant function $t>0$.}
\end{proposition}

\begin{proof}
{
By the Sobolev embedding, $\Phi\in C(H^1(M))$, and $D$ is an open cone. It is nonempty: a nonzero smooth function supported in a small geodesic ball contained in $\{K>0\}$ satisfies $F(u)>0$ and, by the local Poincar\'e inequality, $Q(u)>0$. Let $p=2^*$. For $u\in D$, the fiber $\varphi_u(t)=\frac12Q(u)t^2-\frac1pF(u)t^p$ has a unique critical point $t^-(u)=(Q(u)/F(u))^{1/(p-2)}$, which is a strict global maximizer at the level $\frac1nQ(u)^{n/2}F(u)^{-(n-2)/2}>0$, and $\varphi_u(t)\to-\infty$ as $t\to\infty$. Hence (H1')--(H3') hold in every path-connected component $W$ of $D$.

The set $Y$ enters the proof of Theorem \ref{thmmain} only through the fact that $W$ is an open cone and through \eqref{nomax}, which is used for $v\in\partial W\setminus\{0\}$ and $\varphi_v(\tau)>0$. Since $D$ is an open cone, so is $W$, and $\partial W\cap D=\emptyset$, so it suffices to show that if $v\notin D$, $0<\tau<s$ and $\varphi_v(\tau)>0$, then $\varphi_v(s)\ge\varphi_v(\tau)$. Since $\varphi_v(\tau)>0$, we have $F(v)\le0$: otherwise $Q(v)\le0$ and $\varphi_v<0$. Thus $t\mapsto\varphi_v(t)/t^2=\frac12Q(v)-\frac1pF(v)t^{p-2}$ is nondecreasing and positive at $\tau$, so that $\varphi_v(s)\ge\varphi_v(\tau)$. The proof of Theorem \ref{thmmain} therefore applies to $W$, and yields the minimax identity.

Finally, if $\int_MK\ge0$, then for each constant $t>0$ we have $\Phi(t)=\frac12t^2R_g|M|-\frac1p t^p\int_MK<0$.
}
\end{proof}

\begin{remark}
The last assertion of Proposition \ref{curv} illustrates Remark \ref{local}. If $\int_MK\ge0$, the constants form a ray along which $\Phi<0$. Hence $\Phi$ cannot have the mountain pass geometry with respect to a cone containing $e$ as soon as the direction of the constants can be joined to $e/\|e\|$ by a path in the unit sphere along which $F>0$. This is consistent with the classical necessary condition for \eqref{yamabe}: multiplying \eqref{yamabe} by $u^{1-2^*}$ and integrating, one sees that a positive solution satisfies
\begin{equation*}
	\int_MK=a_n(1-2^*)\int_Mu^{-2^*}|\nabla u|^2+R_g\int_Mu^{2-2^*}<0.
\end{equation*}
To obtain $c_W^->0$, and hence the mountain pass geometry of Theorem \ref{thmMPG}, one needs an analogue of Lemma \ref{li3}, namely a positive lower bound for $Q$ on the unit sphere intersected with $\{F\ge 0\}$. At the critical exponent this is more delicate, since minimizing sequences may concentrate at points where $K>0$. Existence results for \eqref{yamabe} with sign-changing $K$ were obtained by Rauzy \cite{Rau} when $K$ is not too positive; see also Mayer and Zhu \cite{MZ} for a recent variational approach.
\end{remark}

\end{document}